\documentclass[11pt]{article}
\pdfoutput=1 
\usepackage{amsfonts,amssymb,amsmath,amsthm}
\usepackage{algorithm}
\usepackage{algorithmic}
\usepackage{fullpage}
\usepackage{graphicx}
\graphicspath{{./figures/}}  
\usepackage{float}
\usepackage{subfigure}   
\newcommand{\R}{\mathbb{R}}

\providecommand{\He}{\nabla^2 f} 
\newcommand{\eqdef}{\overset{\text{def}}{=}}
\usepackage[colorinlistoftodos,bordercolor=orange,backgroundcolor=orange!20,linecolor=orange,textsize=scriptsize]{todonotes}

\renewcommand{\t}{^{^T}}

\newcommand{\half}{\frac{1}{2}}
\newcommand{\df}{\delta_f}

\providecommand{\Tr}[1]{\mathbf{Tr}\left( #1\right)}

\providecommand{\E}[1]{\mathbb{E}\left[ #1\right]}

\providecommand{\dotprod}[1]{\left< #1\right>}
\providecommand{\norm}[1]{\left\lVert#1\right\rVert}

\newtheorem{lemma}{Lemma}
\newtheorem{theorem}{Theorem}
\newtheorem{assumption}{Assumption}
\usepackage[firstinits=true,backend=bibtex, doi=false,isbn=false,url=false,eprint=false,maxbibnames=10]{biblatex} 
\renewbibmacro{in:}{\ifentrytype{article}{}{\printtext{\bibstring{in}\intitlepunct}}} 
\newcommand{\Ref}[1]{#1}

\DeclareFieldFormat{sentencecase}{\MakeSentenceCase{#1}}

\renewbibmacro*{title}{%
  \ifthenelse{\iffieldundef{title}\AND\iffieldundef{subtitle}}
    {}
    {\ifthenelse{\ifentrytype{article}\OR\ifentrytype{inbook}%
      \OR\ifentrytype{incollection}\OR\ifentrytype{inproceedings}%
      \OR\ifentrytype{inreference}}
      {\printtext[title]{%
        \printfield[sentencecase]{title}%
        \setunit{\subtitlepunct}%
        \printfield[sentencecase]{subtitle}}}%
      {\printtext[title]{%
        \printfield[titlecase]{title}%
        \setunit{\subtitlepunct}%
        \printfield[titlecase]{subtitle}}}%
     \newunit}%
  \printfield{titleaddon}}

\bibliography{\Ref{jabref_library}}

\title{Stochastic Block BFGS: Squeezing More Curvature out of Data}

\author{Robert M. Gower\thanks{School of Mathematics, University of Edinburgh, gower.robert@gmail.com} \qquad Donald Goldfarb\thanks{Department of Industrial Engineering and Operations Research,
Columbia University, New York, goldfarb@columbia.edu} \qquad  Peter Richt\'{a}rik\thanks{School of Mathematics, University of Edinburgh, peter.richtarik@ed.ac.uk. This author would like to acknowledge support from  EPSRC Grant EP/K02325X/1 and  EPSRC Fellowship EP/N005538/1. }}
\begin{document}

\maketitle

\begin{abstract}
We propose a novel limited-memory stochastic block BFGS update for incorporating enriched curvature information in stochastic approximation methods.
In our method, the estimate of the inverse Hessian matrix that is maintained by it, is updated at each iteration using a sketch of the Hessian, i.e., a randomly generated compressed form of the Hessian.
 We propose several sketching strategies, present a new quasi-Newton method that uses stochastic block BFGS updates combined with the variance reduction approach SVRG  to compute batch stochastic gradients, and prove linear convergence of the resulting method. Numerical tests on large-scale logistic regression problems  reveal that our method is more robust and substantially outperforms current state-of-the-art methods.

%
%
%
%
%
%
%
%

\end{abstract}

\section{Introduction}
We design a new stochastic variable-metric (quasi-Newton) method---the stochastic block BFGS method---for solving the Empirical Risk Minimization (ERM) problem:
\begin{equation}\label{eq:prob}
 \min_{x \in \R^d } f(x) \eqdef \frac{1}{n} \sum_{i=1}^n f_i(x).
 \end{equation}
We assume the loss functions $f_i:\R^d\to \R$ to be convex and twice differentiable and focus on the setting where the number of \emph{data points (examples)} ($n$) is very large.

To solve~\eqref{eq:prob}, we employ iterative methods of the form
 \begin{equation}\label{eq:variancereduced}
x_{t+1} = x_{t} - \eta H_{t} g_t,
\end{equation}
where  $\eta >0$ is a stepsize, $g_t\in \R^d$ is an estimate of the gradient $\nabla f(x_t)$ and $H_t \in \R^{d\times d}$ is  a positive definite estimate of the inverse Hessian matrix, that is $H_t \approx \nabla^2 f(x_t)^{-1}$. We refer to $H_t$ as the {\em metric matrix}~\footnote{
Methods of the form~\eqref{eq:variancereduced} can be seen
as estimates of a gradient descent method under the metric defined by $H_t.$ Indeed,
let $\dotprod{x,y}_{H_t} \eqdef \dotprod{H_t^{-1} x,y}$ for any $x,y \in \R^d$ denote an inner product, then the gradient of $f(x_t)$ in this metric is $H_t\nabla f(x_t).$
}.

The most successful classical optimization methods fit the format~\eqref{eq:variancereduced}, such as gradient descent $\left(H_t = I\right)$, Newton's method  $\left(H_t =\nabla^2 f(x_t)^{-1}\right)$, and the quasi-Newton methods $\left(H_t \approx \nabla^2 f(x_t)^{-1}\right)$; all with $g_t =\nabla f(x_t)$. The difficulty
in our setting is that the large number of data points makes the computational costs of a single iteration of these classical methods prohibitively expensive.

%

To amortize these costs, the current state-of-the-art methods use subsampling,  where  $g_t$ and $H_t$ are calculated using only derivatives of the subsampled function
\[
 f_S(x) \eqdef \frac{1}{|S|}\sum_{i\in S}  f_i(x),
\]
where $S \subseteq [n]\eqdef \{1,2,\dots,n\}$ is a subset of examples selected uniformly at random.
Using the subsampled gradient $ \nabla f_S(x)$ as a proxy for the gradient is the basis for
(minibatch) stochastic gradient descent (SGD), but also for  many successful variance-reduced methods~\cite{SAG,SDCA,Johnson2013,S2GD,SAGA,Shalev-Shwartz2013b,MS2GD} that make use of
the subsampled gradient in calculating $g_t.$

Recently,  there has been an effort to  calculate $H_t$ using subsampled Hessian matrices
$\nabla^2 f_T(x_t)$, where $T \subseteq [n]$ is sampled uniformly at random and independently of $S$ \cite{Erdogdu2015,Roosta-Khorasani2016}. The major difficulty in this approach is that calculating $\nabla^2 f_T(x_t)$ can be computationally expensive and, when $d$ is large, storing $\nabla^2 f_T(x_t)$ in memory can be infeasible.
One fairly successful solution to these issues \cite{Byrd2015,Moritz2015} is to use a single Hessian-vector product $\nabla^2 f_T(x_t) v$, where $v \in \R^d$ is a suitably selected vector,  to update $H_t$ by the limited-memory (L-BFGS)~\cite{Nocedal1980} version of the classical BFGS~\cite{Broyden1967,Fletcher1970,Goldfarb1970,Shanno1971} method. Calculating this Hessian-vector product can be done inexpensively using the directional derivative
\begin{equation} \label{eq:Hesssketch}
 \nabla^2 f_T(x_t)v  =\left.\frac{d}{d \alpha}\nabla f_{T}(x_t+\alpha v )\right|_{\alpha=0}.
\end{equation}
In particular, when using automatic differentiation techniques~\cite{Christianson:1992,Griewank2008} or backpropagation on a neural network~\cite{Pearlmutter1994}, evaluating the above costs at most five times as much as the cost of evaluating the subsampled gradient $\nabla f_{T}(x_t)$.


Using only a single Hessian-vector product to update $H_t$ yields only a very limited amount of curvature information, and thus may result in an ineffective metric matrix. The block BFGS method addresses this issue. The starting point in the development of the block BFGS method is the simple observation that, ideally, we would like the metric matrix $H_t$ to satisfy the inverse equation
\[H_t \He_{T}(x_{t}) = I, \]
since then $H_t$ would be the inverse of an unbiased estimate of the Hessian. But solving the inverse equation is computationally expensive. So instead, we propose that $H_t$ should instead satisfy a \emph{sketched} version of this equation, namely
\begin{equation}\label{eq:sketchinver} H_t \He_{T}(x_{t}) D_t= D_t,
\end{equation}
where $D_t \in \R^{d\times q}$ is a randomly generated matrix which has relatively few columns ($q\ll d$). The \emph{sketched subsampled Hessian} $ \He_{T}(x_{t}) D_t$ can be calculated efficiently through $q$ directional derivatives of the form~\eqref{eq:Hesssketch}.

Note that~\eqref{eq:sketchinver} has possibly an infinite number of solutions, including the inverse of $\He_{T}(x_{t})$. To determine $H_t$ uniquely,  we maintain a previous estimate $H_{t-1} \in \R^{d\times d}$ and project $H_{t-1}$ onto the space of symmetric matrices that satisfy~\eqref{eq:sketchinver}. The resulting update, applied to $H_{t-1}$ to arrive at $H_t$, is the block BFGS update.

In the remainder of the paper we describe the block BFGS update, propose a new limited-memory block BFGS update and introduce several new sketching strategies. We  conclude by presenting the results of numerical tests of a method that combines the limited-memory block BFGS update with the SVRG/S2GD method \cite{Johnson2013, S2GD}, and demonstrate that our new method yields dramatically better results when compared to SVRG, or the SVRG method coupled with the classical L-BFGS update as proposed in~\cite{Moritz2015}.

\subsection{Contributions}
This paper makes five main contributions:

\paragraph{\bf (a) New metric learning framework.}
We develop a {\em stochastic block BFGS update}\footnote{In this paper we use the word ``update'' to denote metric learning, i.e., an algorithm for updating one positive definite matrix into another.} for approximately tracking the inverse of the Hessian of $f$. This technique is novel in two ways: the update is more flexible than the traditional BFGS update as it works by employing the actions of a subsampled Hessian on a {\em set} of \emph{random vectors} rather than just on a  single deterministic vector, as is the case with the standard BFGS. That is, we use {\em block sketches} of the subsampled Hessian.

\paragraph{\bf (b) Stochastic block BFGS method.} Our block BFGS update is capable of incorporating {\em enriched} second-order information into gradient based stochastic approximation methods for solving problem \eqref{eq:prob}. In this paper we illustrate the power of this strategy in conjunction with the strategy employed in the SVRG method of \cite{Johnson2013} for computing variance-reduced stochastic gradients. We prove that the resulting combined method is linearly convergent and empirically demonstrate its ability to substantially outperform current state-of-the-art methods.


\paragraph{(c) Limited-memory method.}

To make the stochastic block BFGS method applicable to large-scale problems, we devise a new limited-memory variant of it. As is the case for L-BFGS~\cite{Nocedal1980}, our limited-memory approach allows for a user-defined amount of memory to be set aside. But unlike L-BFGS, our limited-memory approach allows one to use the available memory to store more recent curvature information, encoded in sketches of previous Hessian matrices. This  development of a new limited block BFGS method should also be of general interest to the optimization community.

\paragraph{\bf (d) Factored form.}
We develop  a limited-memory factored form of the block BFGS update.
While factorized versions of standard quasi-Newton methods were proposed in the 1970's (e.g., see~\cite{Murray1972,Goldfarb:1976}), as far as we know, no efficient limited memory versions of such methods have been developed. Factored forms are important as they can be used to enforce positive definiteness of the metric even in the presence numerical imprecision. Furthermore, we use the factored form in calculating new sketching matrices.

\paragraph{\bf (e) Adaptive sketching.}
Not only can sketching be used to tackle the large dimensions of the Hessian, but it can also simultaneously precondition the inverse equation~\eqref{eq:sketchinver}. We present a self-conditioning (i.e., adaptive) sketching that makes use of the efficient factored form of the block BFGS method developed earlier. We also present a sketching approach based on using previous search directions. Our numerical tests show that adaptive sketching can in practice lead to a significant speedup in comparison with sketching from a fixed distribution.

\subsection{Background and Related Work}

The first stochastic variable-metric method developed that makes use of subsampling was the online L-BFGS method~\cite{Schraudolph2006}. In this work the authors adapt the L-BFGS method to make use of subsampled gradients, among other empirically verified improvements. The regularized BFGS method~\cite{Mokhtari2013,Mokhtari2014} also makes use of stochastic gradients, and further modifies the BFGS update by adding a regularizer to the metric matrix.

The first method to use subsampled Hessian-vector products in the BFGS update, as opposed to using differences of stochastic gradients, was the SQN method~\cite{Byrd2015}.
Recently, \cite{Moritz2015} propose combining  SQN with  SVRG. The resulting method performs very well in numerical tests.
In our work we combine a novel {\em stochastic block BFGS} update with SVRG, and prove linear convergence. The resulting method is more versatile and superior in practice to SQN  as it can capture more useful curvature information.

The update formula, that we refer to as the block BFGS update, has a rather interesting background. The formula first appeared in 1983 in unpublished work by Schnabel~\cite{Schnabel1983} on designing quasi-Newton methods that make use of multiple secant equations. Schnabel's method requires several modifications that stem from the lack of symmetry and positive definiteness of the resulting update. Later, and completely independently, the block BFGS update appears in the domain decomposition literature~\cite{Mandel1993} as a preconditioner, where it is referred to as the balancing preconditioner. In that work, the motivation and derivation are very different from those used in the quasi-Newton literature; for instance, no variational interpretation is given for the method. The balancing preconditioner was subsequently taken out of the PDE context and tested as a general purpose preconditioner for solving a single linear system and systems with a changing right hand side~\cite{Gratton2011}. Furthermore, in~\cite{Gratton2011} the authors present a factored form of the update in a different context, which we adapt for limited-memory implementation. Finally, and again independently, a family of block quasi-Newton methods that includes the block BFGS is presented in~\cite{Gower2016} through a variational formulation and in~\cite{Hennig2015} using Bayesian inference.

\section{Stochastic Block BFGS Update}

The stochastic block BFGS update, applied to $H_{t-1}$ to obtain $H_t$, is defined by the projection
\begin{align} H_{t} &= \arg \min_{H \in \R^{d\times d}} \norm{H- H_{t-1}}_{t}^2 \nonumber \\
 &\quad \mbox{subject to }\,\,   H\He_{T}(x_{t})D_{t} =D_{t},  \,\, H = H^T, \label{eq:blockbfgs}
\end{align}
where \[\norm{H}_{t}^2 \eqdef \Tr{H\He_{T}(x_{t})H^T\He_{T}(x_{t})},\]
and $\Tr{\cdot}$ denotes the trace. The method is stochastic since $D_t \in \R^{d \times q}$ is a random matrix. The constraint in~\eqref{eq:blockbfgs} serves as a \emph{fidelity} term, enforcing that $H_{t}$ be a symmetric matrix that satisfies a sketch of the inverse equation~\eqref{eq:sketchinver}. The objective in~\eqref{eq:blockbfgs} acts as a \emph{regularizer}, and ensures that the difference between $H_t$ and $H_{t-1}$ is a low rank update. The choice of the objective is special yet in another sense: recent results show that if the iteration \eqref{eq:blockbfgs} is applied to a fixed invertible matrix $A$ in place of the subsampled Hessian, then the matrices $H_t$ converges to the inverse $A^{-1}$ at a linear rate~\cite{Gower2016}. This is yet one more reason to expect that in our setting the matrices $H_t$ approximately track the inverse Hessian. For similar techniques applied to the problem of solving linear systems, we refer to \cite{Gower2015b, Gower2015c}.

The solution to~\eqref{eq:blockbfgs} is
\begin{eqnarray}
H_{t}  &=& D_{t}\Delta_t D_{t}^T  + \left(I-D_{t}\Delta_t Y_t^T\right) H_{t-1} \left(I-Y_t \Delta_t D_{t}\right),
 \label{eq:BBFGS}
\end{eqnarray}
where $\Delta_t \eqdef (D_t^TY_t)^{-1}$ and $Y_t \eqdef \He_{T}(x_{t})D_{t}.$
This solution was given in~\cite{Gower2014c,Hennig2015} and in~\cite{Schnabel1983} for multiple secant equations.

Note that  stochastic block BFGS~\eqref{eq:BBFGS} yields the same matrix $H_t$ if $D_t$ is replaced by any other matrix $\tilde{D}_t \in \R^{d\times q}$ such that
$\mathbf{span}(\tilde{D}_t) = \mathbf{span}(D_t)$.
It should also be pointed out that the matrix $H_t$ produced by~\eqref{eq:BBFGS} is not what would be generated by a sequence of $q$ rank-two BFGS updates based on the $q$ columns of $D_t$ unless the columns of $D_t$ are $\He_{T}(x_{t})$-conjugate. Note that these columns would be conjugate if they were generated by the BFGS method applied to  the problem of minimizing a strictly convex quadratic function with Hessian $\nabla^2 f_T(x_t)$, using exact line-search.

We take this opportunity to point out that stochastic block BFGS can generate the metric used in the Stochastic Dual Newton Ascent (SDNA) method~\cite{Qu2015} as a special case.
Indeed, when $H_{t-1}=0$, then $H_t$ is given by
\begin{eqnarray} \label{eq:SDNA}
H_t & = & D_{t}(D_{t}^T\He_{T}(x_{t})D_{t})^{-1}D_{t}^T. \end{eqnarray}
When the sketching matrix $D_t$ is a random column submatrix of the identity, then~\eqref{eq:SDNA} is the
positive semidefinite matrix used in calculating the iterates of the SDNA method~\cite{Qu2015}. However, SDNA operates in the dual of~
\eqref{eq:prob}.

\section{Stochastic Block BFGS Method}

The goal of this paper is to design a method that uses a low-variance estimate of the gradient, but also gradually incorporates {\em enriched} curvature information. To this end, we propose to combine the stochastic variance reduced gradient (SVRG) approach~\cite{Johnson2013} with our {\em novel stochastic block BFGS update,} described in the previous section. The resulting method is Algorithm~\ref{alg:ML}.

Algorithm~\ref{alg:ML} has an outer loop in $k$ and an inner loop in $t$.
In the outer loop, the \emph{outer} iterate $w_k \in \R^d$ and the full gradient $\nabla f(w_k)$ are computed. In the inner loop, both the estimate of the gradient $g_t$ and our metric $H_t$ are updated using the SVRG update and the block BFGS update, respectively.

To form the sketching matrix $D_{t}$ we employ one of the following three strategies:

\paragraph{(a) Gaussian sketch.} $D_{t}$ has standard Gaussian entries sampled i.i.d at each iteration.

\paragraph{(b) Previous search directions delayed.} Let us write $d_t=-H_tg_t$ for the search direction used in step $t$ of the method. In this strategy we store $L$ such search directions as columns of matrix $D_t$: $D_{t} = [d_{t+1-L},\ldots, d_t]$, and then update $H_{t}$ only once every $L$ inner iterations.

\paragraph{(c) Self-conditioning.}  We sample $C_t \subseteq [d] $ uniformly at random and set $D_{t} = L_{t-1}I_{:C_t} =[L_{t-1}]_{:C_t}$, where $L_{t-1}L_{t-1}^T=H_{t-1}$ and $I_{:C_t}$ denotes the concatenation of the columns of the identity matrix indexed by a set $C_t \subset [d]$. Thus the sketching matrix is formed with a random subset of columns of a factored form of $H_t.$ The idea behind this strategy is that the ideal sketching matrix should be $D_t = (\nabla^2 f_T(x_t))^{-1/2}$ so that the sketch not only compresses but also acts as a preconditioner on the inverse equation~\eqref{eq:sketchinver}. It was also shown in~\cite{Gower2016} that this choice of sketching matrix can accelerate the convergence of $H_t$ to the inverse of a fixed matrix. In Section~\ref{sec:Factored} we describe in detail how to efficiently maintain and update the factored form.

\begin{algorithm}[!ht]
\begin{algorithmic}
\STATE {\bf inputs:} $w_0 \in \R^d$, stepsize $\eta >0$, $s = $ subsample size, $q = $ sample action size, $m = $ size of the inner loop.
 \STATE {\bf initiate:}  $H_{-1} =I$ ($d\times d$ identity matrix)
\FOR {$k = 0, 1, 2, \dots$}
	\STATE Compute the full gradient $\mu = \nabla f(w_{k})$  
		\STATE Set $x_0 = w_k$
	\FOR{$t = 0,\ldots, m-1$}
			   \STATE Sample $S_t,T_t \subseteq [n]$, independently
               \STATE Compute a variance-reduced stochastic gradient:  $g_t = \nabla f_{S_t}(x_t)-\nabla f_{S_t}(x_0) +\mu$
			   \STATE Form $D_t \in \R^{d\times q}$ so that rank$(D_t) = q$
			   \STATE Compute $Y_t = \nabla^2 f_{T_t}(x_t)D_t$ (without ever forming the $d\times d$ matrix $\nabla^2 f_{T_t}(x_t)$) 
	           \STATE Compute ${D_t}^T Y_t$  and its Cholesky factorization 	(this implicitly forms $\Delta_t = (D_t^T Y_t)^{-1}$)	
			    \STATE {\em Option I:} Use \eqref{eq:BBFGS} to obtain $H_t$ and set $d_t = -H_t g_t$
			   \STATE {\em Option II:} Compute $d_t = -H_tg_t$ via Algorithm~\ref{alg:L-BBFGS}
			   \STATE Set $x_{t+1} = x_{t} + \eta d_t$
	\ENDFOR	
	\STATE {\em Option I:} Set $w_{k+1} = x_m$	
	\STATE {\em Option II:} Set $w_{k+1} = x_i$, where $i$ is selected uniformly at random 	from $[m]=\{1,2,\dots,m\}$		
\ENDFOR
 \STATE {\bf output:} $w_{k+1}$
\end{algorithmic}
\caption{Stochastic Block BFGS Method}
\label{alg:ML}
\end{algorithm}

\subsection{Limited-Memory Block BFGS}

When $d$ is large, we cannot store the  $d\times d$ matrix $H_{t}$. Instead, we store $M$ block triples, consisting of previous block curvature pairs and the inverse of their products
\begin{equation}\label{eq:actionpairs}
\left( D_{t+1-M}, Y_{t+1-M}, \Delta_{t+1-M}\right),\ldots, \left( D_{t}, Y_t , \Delta_t \right).\end{equation}
With these triples we can form the $H_{t}$ operator implicitly by using a block limited-memory two loop recurrence. To describe this two loop recurrence, let
$V_{t} \eqdef  I-D_{t}\Delta_t Y_t^T.$

The block BFGS update~\eqref{eq:BBFGS} with memory parameter $M$ can be expanded as a function of the $M$ block triples \eqref{eq:actionpairs} and of $H_{t-M}$ as
\begin{eqnarray*}
H_{t}
&=&  V_{t} H_{t-1}V_{t}^T  +D_{t}\Delta_t D_{t}^T\\
&=& V_{t}\cdots V_{t+1-M} H_{t-M} V_{t+1-M}^T \cdots V_{t}^T + \sum_{i=t}^{t+1-M} V_{t}\cdots V_{i+1} D_{i}\Delta_i D_{i}^T V_{i+1}^T\cdots V_{t}^T.
\end{eqnarray*}
Since we do not store $H_t$ for any $t$, we do not have access to $H_{t-M}$. In our  experiments we simply set $H_{t-M}=I$ (the identity matrix). Other, more sophisticated, choices are possible, but we do not explore them further here. Using the above expansion, the action of the operator $H_{t}$ on a vector $v$ can be efficiently calculated using Algorithm~\ref{alg:L-BBFGS}.

\begin{algorithm}
\begin{algorithmic}
\STATE {\bf inputs:}
$ g_t \in \R^d, D_{i}$, $Y_i \in \R^{d\times q}$ and $\Delta_i \in \R^{q \times q}$
 for $i \in \{t+1-M,\ldots,t\} $.
 \STATE {\bf initiate:} $v = g_t$
\FOR{$i = t, \ldots, t-M+1$}
	\STATE $\alpha_i =  \Delta_i D_i^T v$
	\STATE $v = v-  Y_i \alpha_i$
\ENDFOR
\FOR{$i = t-M+1, \ldots, t$}
	\STATE $\beta_i=  \Delta_iY_i^T v$
	\STATE $v = v + D_i( \alpha_i - \beta_i)$
\ENDFOR
\STATE {\bf output} $H_t g_t = v$
\end{algorithmic}
\caption{Block L-BFGS Update (Two-loop Recursion)}
\label{alg:L-BBFGS}
\end{algorithm}

The total cost  in  floating point operations of executing Algorithm~\ref{alg:L-BBFGS} is $Mq(4d+2q).$ In our experiments $M=5$ and $q$ will be orders of magnitude less than $d$, typically $q \leq \sqrt{d}$. Thus  the cost of applying Algorithm~\ref{alg:L-BBFGS} is approximately $O(d^{3/2}).$ This does not include the cost of computing the product ${D_t}^T Y_t$  ($O(q^2 d)$ operations) and its Cholesky factorization ($O(q^3)$ operations), which is done outside of Algorithm \ref{alg:L-BBFGS}. The two places in Algorithm~\ref{alg:L-BBFGS} where multiplication by $\Delta_i$ is indicated is in practice performed by solving two triangular systems using the Cholesky factor of ${D_i}^T Y_i$. We do this because it is more numerically stable than explicitly calculating the inverse matrix $\Delta_i.$

%


\subsection{Factored Form} \label{sec:Factored}
Here we develop a new efficient method for maintaining and updating a {\em factored form of the metric matrix}. This facilitates the development of  a novel idea which we call  {\em self-conditioning sketch.}

Let $L_{t-1}\in \R^{d\times d}$ be invertible such that $L_{t-1}L_{t-1}^T =H_{t-1}$. Further, let $ G_{t} =(D^T_{t}L_{t-1}^{-T} L_{t-1}^{-1}D_{t})^{1/2} $ and $R_t = \Delta_{t}^{1/2}$. An update formula for the factored form of $H_{t}$, i.e., for $L_{t}$ for which $H_{t} = L_{t} L_{t}^T$, was recently given (in a different context) in~\cite{Gratton2011}:
\begin{equation}\label{eq:Chol}L_{t} =V_{t} L_{t-1}+D_{t}R_{t}G_{t}^{-1}D_{t}^T L_{t-1}^{-T}.
\end{equation}
This factored form of $H_{t}$ is too costly to compute because it requires inverting $L_{t-1}$. However, if we let $D_{t} = L_{t-1}I_{:C_{t}}$, where $C_t \subset [d]$, then~\eqref{eq:Chol} reduces to
\begin{equation}\label{eq:CholLI}
 L_{t}  = V_{t}L_{t-1} +D_{t}R_{t}I_{C_t:}
\end{equation}
which can be computed efficiently. Furthermore, this update of the factored form~\eqref{eq:CholLI} is amenable to recursion and can thus be expanded as
\begin{eqnarray}
 L_{t}  &=& V_{t}\left(V_{t-1}L_{t-2} +D_{t-1}R_{t-1}I_{C_{t-1}:}\right) +D_{t}R_{t}I_{C_t:} \nonumber \\
 & =& V_{t}\cdots V_{t+2-M} V_{t+1-M} L_{t-M} + V_{t}\cdots V_{t+2-M} D_{t+1-M}R_{t+1-M}I_{C_{t+1-M}:}
  + D_{t}R_{t}I_{C_t:} \quad.\label{eq:factexpand}
\end{eqnarray}
By storing $M$ previous curvature pairs~\eqref{eq:actionpairs} and additionally the sets $C_{t+1-M},\ldots, C_t$, we can calculate the action of $L_{t}$ on a matrix $V \in \R^{d\times q}$ by using~\eqref{eq:factexpand}, see Algorithm~\ref{alg:facL-BBFGS}. To the best of our knowledge, this is the first limited-memory factored form in the literature.
Since we do not store $L_t$ explicitly for any $t$ we do not have access to $L_{t-M}$, required in computing~\eqref{eq:factexpand}. Thus we simply use $L_{t-M} =I$.

\begin{algorithm}
\begin{algorithmic}
\STATE {\bf inputs:} $V$, $D_{i},Y_i, \Delta_i \in \R^{d\times q}$ and $C_i \subset [d]$, for $i \in \{t+1-M,\ldots,t\} $.
 \STATE {\bf initiate:} $W = V$ 
\FOR{$i = t+1-M, \ldots, t$}
 \STATE $W  = W -   D_i \Delta_i Y_i^T W +D_iR_iW_{C_i:}$
\ENDFOR
\STATE {\bf output} $W$ (we will have $W = L_tV $)
\end{algorithmic}
\caption{Block L-BFGS Update (Factored loop  recursion for computing $L_tV$)}
\label{alg:facL-BBFGS}
\end{algorithm}
Again, we can implement a more numerically stable  version of Algorithm~\ref{alg:facL-BBFGS} by storing the Cholesky factor of $D_i^TY_i $ and using triangular solves, as opposed to calculating the inverse matrix $\Delta_i =(D_i^TY_i)^{-1} $.


\section{Convergence}

In this section we prove that Algorithm~\ref{alg:ML} converges linearly. Our analysis relies on the following assumption, and is a combination of novel insights and techniques from \cite{S2GD} and \cite{Moritz2015}.
\begin{assumption} \label{ass:strongsmooth}
There exist constants $ 0< \lambda \leq \Lambda $ such that
\begin{equation}\label{eq:hessbound}
\lambda I   \preceq \nabla^2 f_T(x) \preceq  \Lambda I
\end{equation}
for all $x\in \R^d$ and all $T \subseteq [n]$.
\end{assumption}

We need two technical lemmas, whose proofs are given in Sections 7 and 8 at the end of the paper.
\begin{lemma} \label{lem:Hspectra}
Let $H_t$ be the result of applying  the limited-memory Block BFGS update with memory $M$, as implicitly defined by Algorithm~\ref{alg:L-BBFGS}.
Then there exists positive constants $\Gamma \geq \gamma >0$ such that for all $t$ we have
\begin{equation}\label{eq:Hspectra}
\gamma I\preceq H_{t} \preceq \Gamma I.
\end{equation}
\end{lemma}
A proof of this lemma is given in Appendix 2, where in particular it is shown that the lower bound satisfies $\gamma~\geq~\frac{1}{ 1 + M \Lambda}$  and the upper bound satisfies
\begin{equation}\label{eq:Gammakappbnd}
\Gamma~\leq~(1~+~\sqrt{\kappa})^{2M} \left(1 + \frac{1}{\lambda(2\sqrt{\kappa} + \kappa)}\right),
\end{equation} where
$\kappa~\eqdef~\Lambda/\lambda.$

We now state a bound on the norm of the SVRG  variance-reduced gradient for minibatches.
\begin{lemma}\label{lem:pretheo}Suppose Assumption~\ref{ass:strongsmooth} holds, let $w_*$ be the unique minimizer of $f$ and let $w,x \in \R^d.$ Let $\mu =\nabla f(w)$ and  $g = \nabla f_S(x)-\nabla f_S(w)+\mu$. Taking expectation with respect to $S$, we have
\begin{align} \label{eq:lemma2}
\E{\norm{g}_2^2 } \leq 4\Lambda( f(x)-f(w_*)) + 4(\Lambda-\lambda)(f(w)-f(w_*)).
\end{align}
\end{lemma}

The following theorem guarantees the linear convergence of Algorithm~\ref{alg:ML}.
\begin{theorem}
Suppose that Assumption~\ref{ass:strongsmooth} holds. Let $w_*$ be the unique minimizer of $f$. When Option~II is used in Algorithm~\ref{alg:ML}, we have for all $k \geq 0$ that
\[\E{f(w_{k})-f(w_*)} \leq \rho^k \E{f(w_0)-f(w_*)},\]
where the convergence rate is given by
\[\rho  = \frac{1/2m\eta+\eta\Gamma^2 \Lambda(\Lambda-\lambda)}{ \gamma\lambda -\eta\Gamma^2 \Lambda^2}  < 1,\]
assuming we have chosen $\eta < \gamma\lambda/(2\Gamma^2\Lambda^2) $ and that we choose $m$ large enough to satisfy\footnote{By our assumption on $\eta$, the expression on the right is nonnegative.}
\[m \geq \frac{1}{2\eta \left(\gamma\lambda- \eta \Gamma^2 \Lambda(2\Lambda-\lambda)\right)}.\]
\end{theorem}
\noindent \emph{Proof.}
 Since $g_{t} = \nabla f_{S_t}(x_t)-\nabla f_{S_t}(w_{k}) +\mu$ and $x_{t+1}=x_t-\eta H_t g_t$ in Algorithm~\ref{alg:ML}, from~\eqref{eq:hessbound} we have that
\begin{eqnarray*}
f(x_{t+1})
&\leq &
f(x_t) + \eta \nabla f(x_t)^Td_t+\frac{\eta^2\Lambda}{2}\norm{d_t}_{2}^2\\
&= &
f(x_t)  -\eta \nabla f(x_t)^T H_t g_{t}+\frac{\eta^2\Lambda}{2}\norm{H_t g_{t}}_{2}^2.
\end{eqnarray*}
Taking expectation conditioned on $x_t$ (i.e., with respect to $S_t$, $T_t$ and $D_t$)  and using Lemma~\ref{lem:Hspectra} we have
\begin{eqnarray}
\E{f(x_{t+1}) \, | \, x_t} &\leq &  f(x_t)  -\eta \E{\nabla f(x_t)^T H_t \nabla f(x_t)\, | \, x_t}+\frac{\eta^2\Lambda}{2}\E{\norm{H_t g_{t}}_{2}^2\, | \, x_t} \nonumber \\
&\overset{\text{Lemma~\ref{lem:Hspectra}}}{\leq} &
   f(x_t)  -\eta \gamma \norm{\nabla f(x_t)}_2^2+\frac{\eta^2\Gamma^2 \Lambda}{2}\E{\norm{ g_{t}}_2^2 \, | \, x_t}. \label{eq:theorproof1}
\end{eqnarray}
Introducing the notation $\df(x) \eqdef f(x) - f( w_*)$  and applying Lemma~\ref{lem:pretheo} and the fact that strongly convex functions satisfy the inequality
$\|\nabla f(x)\|_2^2 \geq 2\lambda \df(x)$ for all $x\in \R^d$,  gives
\begin{eqnarray*}
\E{f(x_{t+1}) \, | \, x_t}
 &\overset{\eqref{eq:theorproof1}+\eqref{eq:lemma2}}{\leq} &
  f(x_t) -2\eta \gamma\lambda \df(x_t) +2\eta^2\Gamma^2 \Lambda\left(\Lambda
\df(x_t) ) + (\Lambda-\lambda) \df(w_k) )\right) \\
& = & f(x_t) - \alpha \df(x_t) + \beta \df(w_k).
\end{eqnarray*}
where $\alpha = 2\eta \left(\gamma\lambda -\eta\Gamma^2 \Lambda^2\right)$ and  $\beta = 2\eta^2\Gamma^2 \Lambda(\Lambda-\lambda)$.
Taking expectation, summing over $t =0, \ldots, m-1$ and using telescopic cancellation gives
\begin{eqnarray*}
 \E{f(x_{m})}   &= & \E{f(x_0)} -\alpha \left(\sum_{t=1}^{m-1}\E{\df(x_t)} \right) +m\beta \E{\df(w_k)}\\
&=& \E{f(w_{k})}
-m \alpha\E{\df(w_{k+1})} + m\beta \E{\df(w_k)},
\end{eqnarray*}
where we used that $w_{k} =x_0$ and  $\sum_{t=1}^m  \E{x_t} = m \E{w_{k+1}}$ which is a consequence of using  Option II in Algorithm~\ref{alg:ML}.  Rearranging the above gives
\begin{eqnarray*}
0 &\leq & \E{f(w_{k}) -f(x_m)}  -m\alpha \E{\df(w_{k+1})} +m\beta \E{\df(w_k)}\\
&\leq & \E{\df(w_k)} -m\alpha \E{\df(w_{k+1})} +m\beta \E{\df(w_k)}\\
&= & -m\alpha\E{\df(w_{k+1})} +(1+m\beta) \E{\df(w_k)}
\end{eqnarray*}
where we used that $f(w_*) \leq f(x_m)$. Using that $\eta~<~\gamma \lambda/(2\Gamma^2\Lambda^2),$ it follows that
\[
\E{\df(w_{k+1})}
\leq \frac{1+2m\eta^2\Gamma^2 \Lambda(\Lambda-\lambda)}{2m\eta \left( \gamma\lambda -\eta\Gamma^2 \Lambda^2\right)}\E{\df(w_k)}. \quad \qed
\]

\section{Numerical Experiments}

To validate our approach, we compared our algorithm to SVRG~\cite{Johnson2013}  and the variable-metric adaption of SVRG presented in~\cite{Moritz2015},  which we refer to as the \texttt{MNJ} method. We tested the methods on seven empirical risk minimization problems with a logistic loss and L2 regularizer, that is, we solved
\begin{equation} \label{eq:logprob}
 \min_{w} \sum_{i=1}^n \ln\left( 1 + \exp(-y_i \dotprod{a^i,w})\right) + \frac{1}{n} \norm{w}_2^2,
\end{equation}
where $A =[a^1, \ldots, a^n] \in \R^{d\times n}$ and $y \in \{0, 1 \}^n$ are the given data. We have also employed the standard ``bias trick''\footnote{The bias trick is to add an additional \emph{bias} variable $\beta\in \R$ so that the exponent in~\eqref{eq:logprob} is  $-y_i (\dotprod{a^i,w} + \beta)$. This is done efficiently by simple concatenating a row of ones to data matrix so that $\dotprod{ [a^i \, 1], [w \, \beta]} = \dotprod{a^i,w} + \beta,$ for $i=1,\ldots, n.$}.
For our experiments we used data from the LIBSVM collection~\cite{Chang2011}.
All the methods were implemented in MATLAB. All the code for the experiments can be downloaded from \url{http://www.maths.ed.ac.uk/~prichtar/i_software.html}.


We tested three variants of Algorithm~\ref{alg:ML}, each specified by the use of  a different sketching matrix. 
In Table~\ref{tab:labelkey} we present a key to the abbreviations used in all our figures. The first three methods, \texttt{gauss}, \texttt{prev} and \texttt{fact}, are 
implementations of the three variants (a), (b) and (c), respectively,  of Algorithm~\ref{alg:ML} using the three different sketching methods discussed at the start of Section 3. In these three methods the $M$ stands for the number of stored curvature pairs~\eqref{eq:actionpairs} used.

For each method and a given parameter choice, we tried the stepsizes
 \[\eta \in \{10^0,\; 5\cdot 10^{-1},\; 10^{-1},\;\ldots ,\; 10^{-7}, \;5\cdot 10^{-8}, \;10^{-8} \}\] 
 and reported the one that gave the best results.
 We used the error $f(x_t)-f(w_*)$ for the $y$-axis in all our figures\footnote{Thanks to Mark Schmidt, whose code \texttt{prettyPlot} 
  was  used to generate all figures: \url{https://www.cs.ubc.ca/~schmidtm/Software/prettyPlot.html}. Note that in all plots the markers (circles, plus signs, triangles ... etc) are used to help distinguish the different methods, and are not related to the iteration count.}.
  We calculated   $f(w_*)$ by running all methods for 30 passes over the data, and then taking the minimum function value.

\begin{table} \centering
\begin{tabular}{| c | c |}
\hline
Method & Description \\[0.25cm]
\hline
\texttt{gauss}\_$q$\_$M$ & \parbox{7cm}{$D_t\in \R^{d \times q}$ with i.i.d Gaussian entries}\\[0.25cm]
\texttt{prev}\_$L$\_$M$ &  \parbox{7cm}{$D_t = [d_t, \ldots, d_{t-L+1}]$. Updated every $L$ inner iterations}\\[0.25cm]
\texttt{fact}\_$q$\_$M$ &  \parbox{7cm}{$D_t = L_{t-1}I_{:C}$ where $C \subset [n]$  sampled uniformly at random and $|C|=q$}\\[0.25cm]
\texttt{MNJ}\_$|T_t|$ & \parbox{7cm}{Algorithm~1 in~\cite{Moritz2015} where $|T_t|=$ size of Hessian subsampling and $L =10$}\\[0.25cm]
\hline
\end{tabular}
\caption{A key to the abbreviations used for each method}
\label{tab:labelkey}
\end{table}

Finally, we used $m = \lfloor n/|S_t| \rfloor$  for the number of inner iterations in 
all variants of Algorithm~\ref{alg:ML}, SVRG and \texttt{MNJ}, so that 
all methods perform an entire pass over the data before recalculating the gradient.

\subsection{Parameter investigation}\label{sec:numparam}
\newcommand\widshrink{0.46}
\newcommand\heightshrink{0.40}

In our first set of tests we explored the parameter space of the \texttt{prev} variant of Algorithm 1.
By fixing all parameters but one, we can see how sensitive the \texttt{prev} method is to the singled-out parameter, but also, build some intuition as to what value or, interval of values, yields the best results for this parameter. We focused our tests on the \texttt{prev} method since it proved to be 
overall, the 
most robust method.

\begin{figure}
\centering
\subfigure[memory $M$]{
	\includegraphics[width =  \columnwidth, height =\heightshrink\columnwidth ]{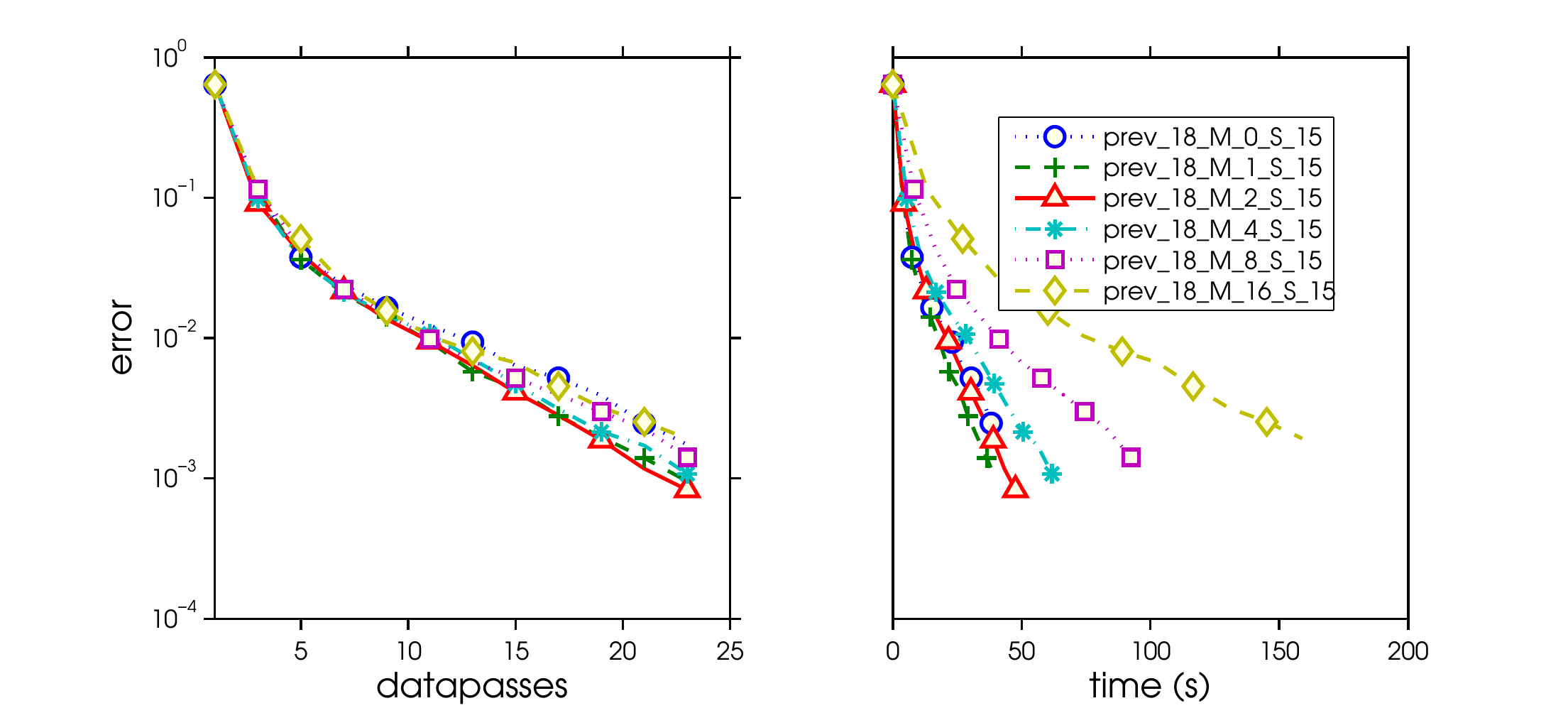}
}\\
\subfigure[Hessian subsampling $|T_t|$]{	
\includegraphics[width =  \columnwidth, height =\heightshrink\columnwidth ]{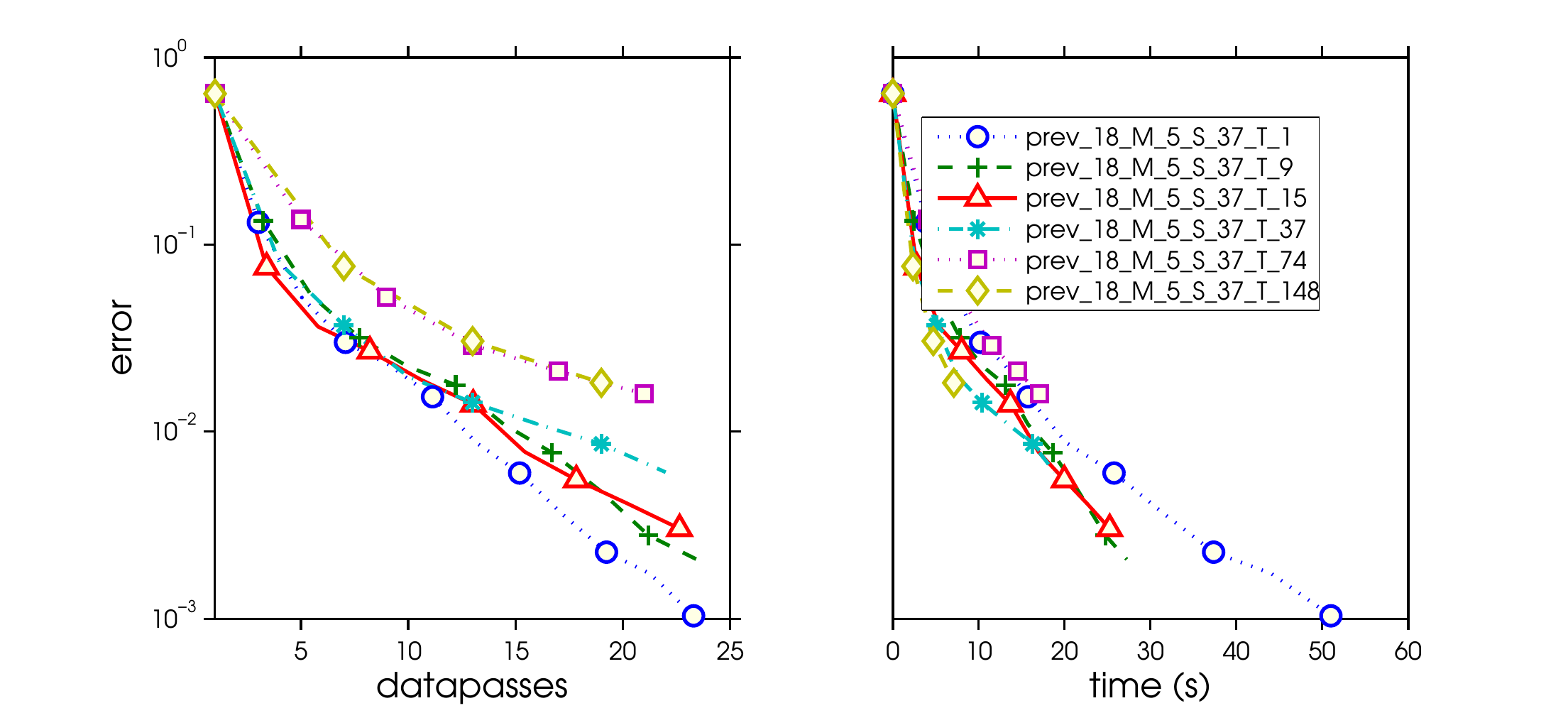}
	}
    \caption{\texttt{w8a}   $(d;n) = (300; \; 49,749)$ }
     \label{fig:param}
\end{figure}

In Figure~\ref{fig:param} (a) we depict the results of varying the memory parameter $M$, while fixing the remaining parameters. In particular, we fixed 
$|T_t|=|S_t| = 15.$ In both the subplots of  error $\times$ datapasses and error $\times$ time in Figure~\ref{fig:param} (a), we see that $1 \leq M\leq 4$ resulted in the best performance. Furthermore, the error $\times$ datapasses is insensitive to increasing the memory, since increasing the memory parameter does not incur in any additional data passes. On larger dimensional problems we found that approximately $M=5$ 
yielded 
the overall best performance. Thus we used $M=5$ in our tests on large-scale  problems in Sections~\ref{sec:numdatapass} and~\ref{sec:numtime}.

In Figure~\ref{fig:param} (b) we 
experimented 
with varying $|T_t|$, the size of the Hessian subsampling. Here the results of the  error $\times$ datapasses subplot 
conflict 
with those of the the error $\times$ time subplots. While in the  error $\times$ time subplot the method improves as $|T_t|$ increases, in the  error $\times$ datapasses subplot the method improves as $|T_t|$ decreases.
As a compromise, in Sections~\ref{sec:numdatapass} and~\ref{sec:numtime} we use
$|S_t| =|T_t|$ as our default choice.

\begin{figure}
\centering
\subfigure[subsampling $|S_t| =|T_t|$]{
	\includegraphics[width =  \columnwidth, height =\heightshrink\columnwidth ]{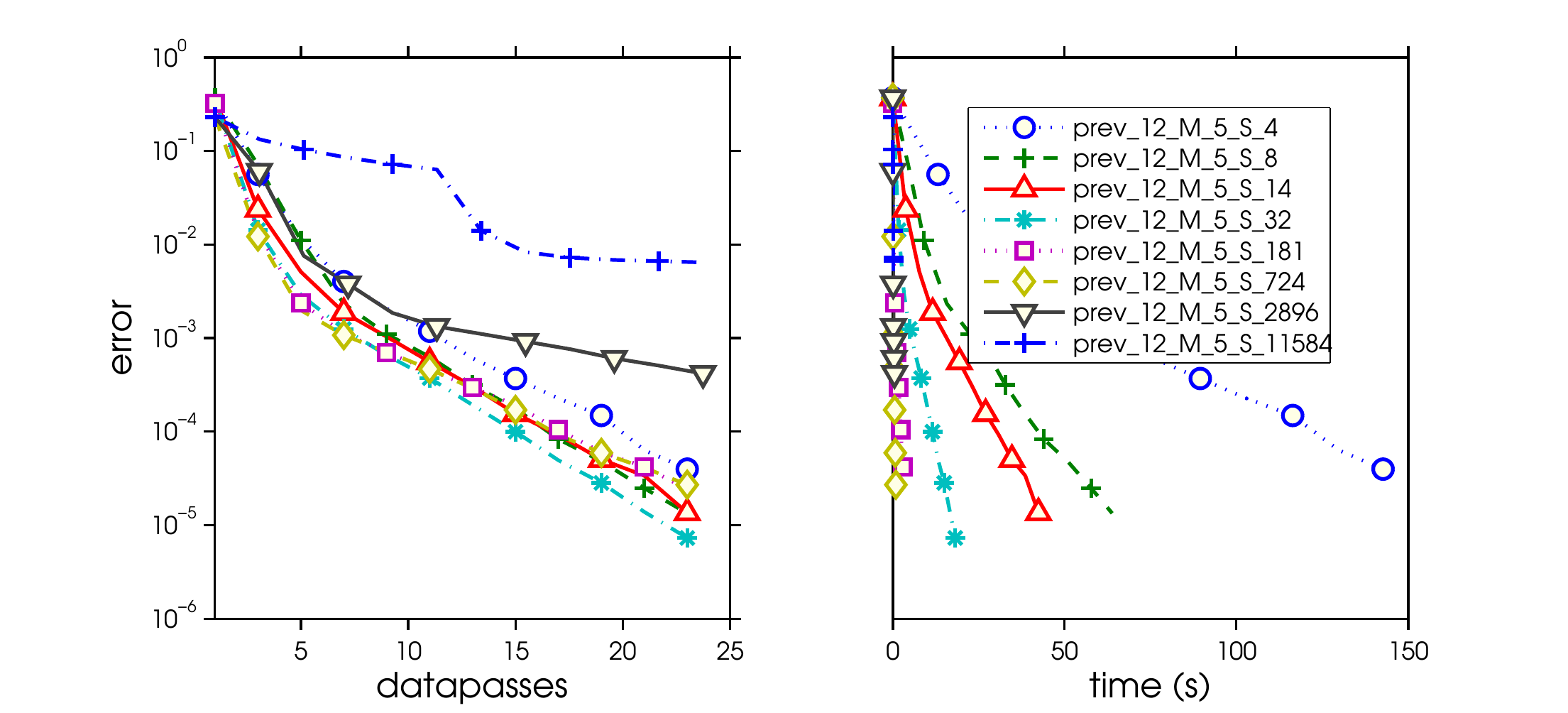}
}\\
\subfigure[update size $L$]{
	\includegraphics[width =  \columnwidth, height =\heightshrink\columnwidth ]{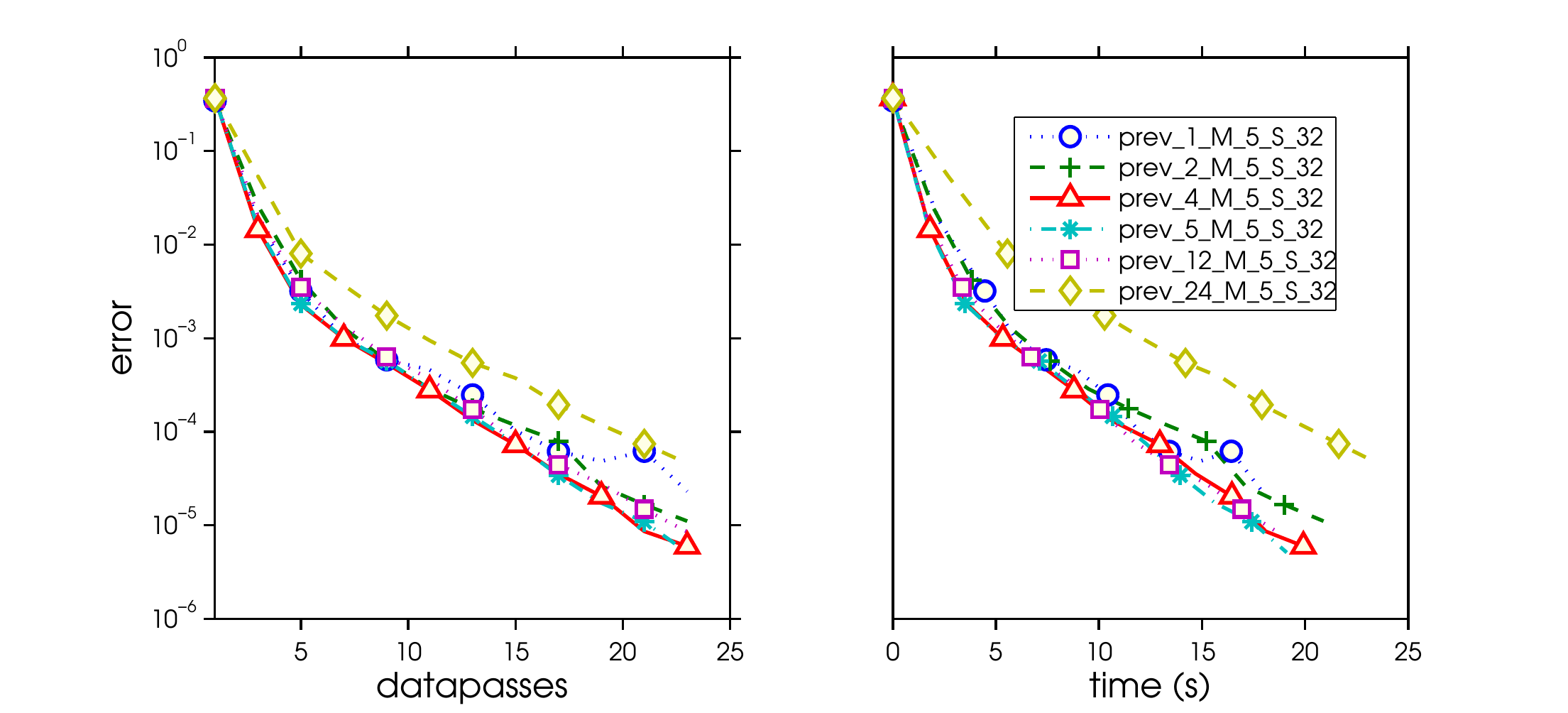}
}%
    \caption{\texttt{a9a}  $(d;n) = (123; \; 32,561)$ } \label{fig:param2}
\end{figure}

In Figure~\ref{fig:param2} (a) we 
experimented 
with varying the size of gradient subsampling and Hessian subsampling jointly with $|S_t|=|T_t|.$  In both the  error $\times$ datapasses  and the error $\times$ time subplot the range
\[\lceil n^{1/3} \rceil =32 \leq |S_t| = |T_t| \leq 724 = \lceil n^{1/2} \rceil,\]
 resulted in a good performance. Based on this experiment, we use $|S_t|=|T_t|= n^{1/2}$ as our default choice in Sections~\ref{sec:numdatapass} and~\ref{sec:numtime}.
 Note that when $|S_t|$ is large the method passes through the data in fewer iterations and consequentially less time. This is why the method terminated early (in time) when $|S_t|$ is large.

Finally in Figure~\ref{fig:param2} (b) we vary the parameter $L$, that is, the number of previous search directions used to form the columns of $D_t$.
From these tests we can see that using a value that is too small, e.g. $L =1$, or a value that is too large $L = 2 \lceil \sqrt{d} \rceil =24$,  results in an inefficient method in terms of both time and datapasses. Instead, we get the best performance when $\lceil d^{1/4} \rceil \leq L \leq  \lceil d^{1/2} \rceil.$
Thus on the first five problems we used either $L = d^{1/3}$ or $L=d^{1/4}$, depending on which gave the best performance. As for the last two problems:  \texttt{rcv1-train.binary} and \texttt{url-combined}, we found that  $L=d^{1/4}$ was too large. Instead, we probed the set $L \in \{2, \ldots, 10 \}$ for an $L$ that resulted in a reasonable performance.

Through these four experiments in Figures~\ref{fig:param} and~\ref{fig:param2},
we can conclude that the \texttt{prev} method is not overly sensitive to the choice of these parameters. That is, the method works well for a range of parameter choices. This is in contrast with choosing the stepsize parameter,  whose %
choice 
can make the difference between a divergent method and a fast method.


\subsection{Data passes}\label{sec:numdatapass}
We now compare all the methods in Table~\ref{tab:labelkey} in terms of error $\times$ datapasses. Since our experiment in Figure~\ref{fig:param} (b) indicated that $|S_t|=|T_t|$ resulted in a reasonable method, for simplicity, we used the same subsampling for the gradient and Hessian in all of our methods; that is $S_t = T_t.$ This is not necessarily an optimal choice.
Furthermore,  we 
set the subsampling size $|S_t| = \sqrt{n}$.

For the \texttt{MNJ} method, we used the suggestion of the authors of both~\cite{Byrd2015} and~\cite{Moritz2015}, and chose $L=10$ and $|T_t|\approx L|S_t|$ so that the computational workload performed by $L$ inner iterations of the SVRG method was approximately equal to that of applying the L-BFGS metric once. The exact rule we found to be efficient was
\[ |T_t| = \left\lfloor \min\left\{ \frac{L  |S_t|}{2}, n^{2/3} \right\} \right \rfloor. \]
We set the memory to $10$ for the \texttt{MNJ} method in all tests, which is a standard choice for the L-BFGS method.


On the problems with $d$ significantly smaller then $n$, such as Figures 1(c) and 1(d), all the methods that make use of curvature information performed similarly and significantly better than the SVRG method.
The \texttt{prev} method proved to be the best overall and the most robust, performing comparably well on problems with $d \ll n$, such as  Figures 1(c) and 1(d), but was also the most efficient method in Figures 1(a), 1(b), 1(e) (shared being 
most efficient with \texttt{MNJ}) and Figure~2(a). The only problem on which the \texttt{prev} method was not the most efficient method was on the \texttt{url-combined} problem in Figure 2(b), where the \texttt{MNJ} method proved to be the most efficient.

Overall, these experiments illustrate that incorporating  curvature information results in a fast and 
robust method.
Moreover, the added flexibility of the block BFGS update to incorporate more curvature information, as compared to a single Hessian-vector product in the \texttt{MNJ} method, can significantly improve the convergence of the method, as can be seen in Figures~1(a) and~1(b).
\begin{figure}
 \centering
    \subfigure[ \texttt{gissette\_scale} ]{
    \includegraphics[width =  \widshrink\columnwidth, height =\heightshrink\columnwidth ]{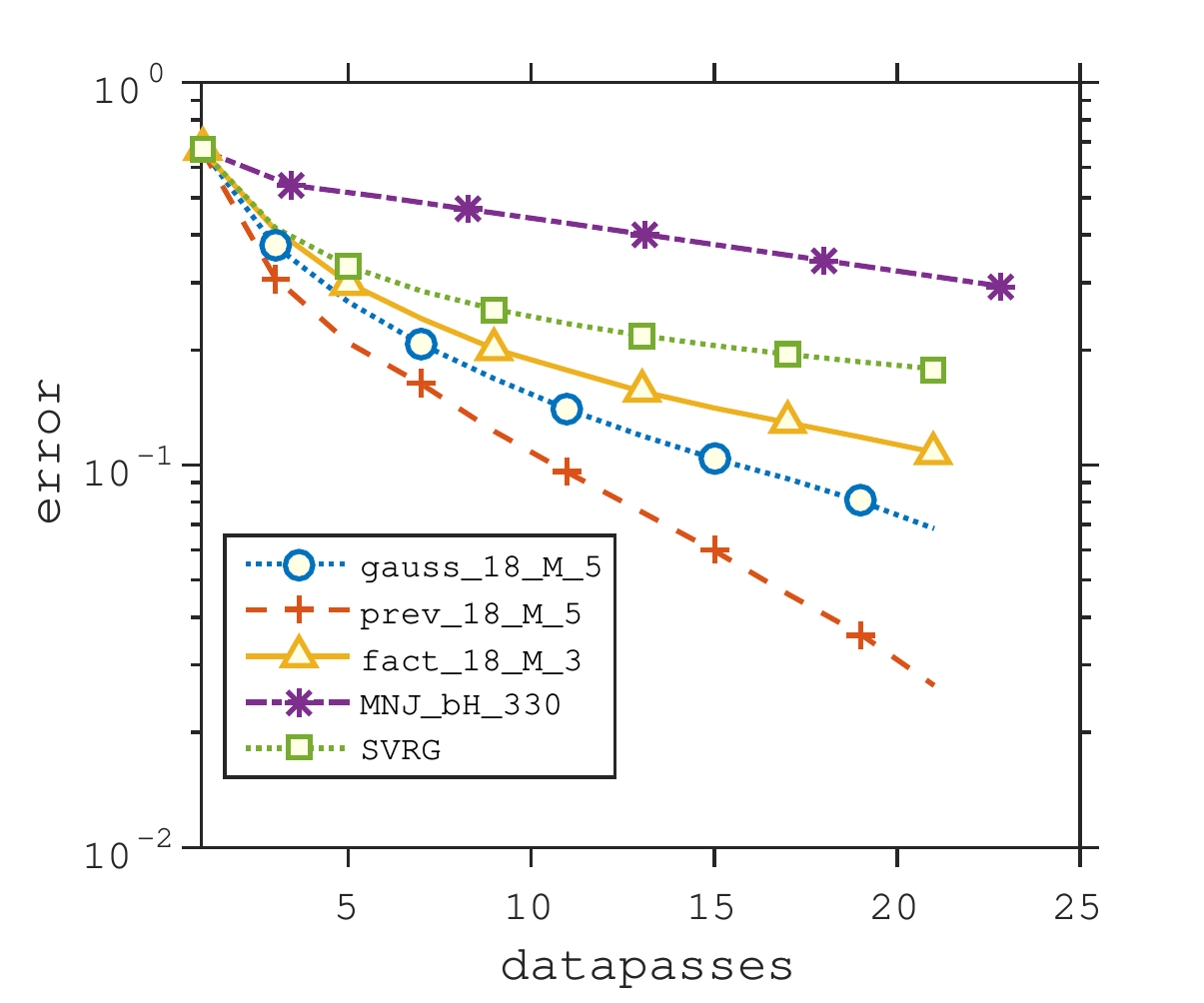}
    }\hfill
    \subfigure[ \texttt{covtype-libsvm-binary} ]{
    \includegraphics[width =  \widshrink\columnwidth, height =\heightshrink\columnwidth ]{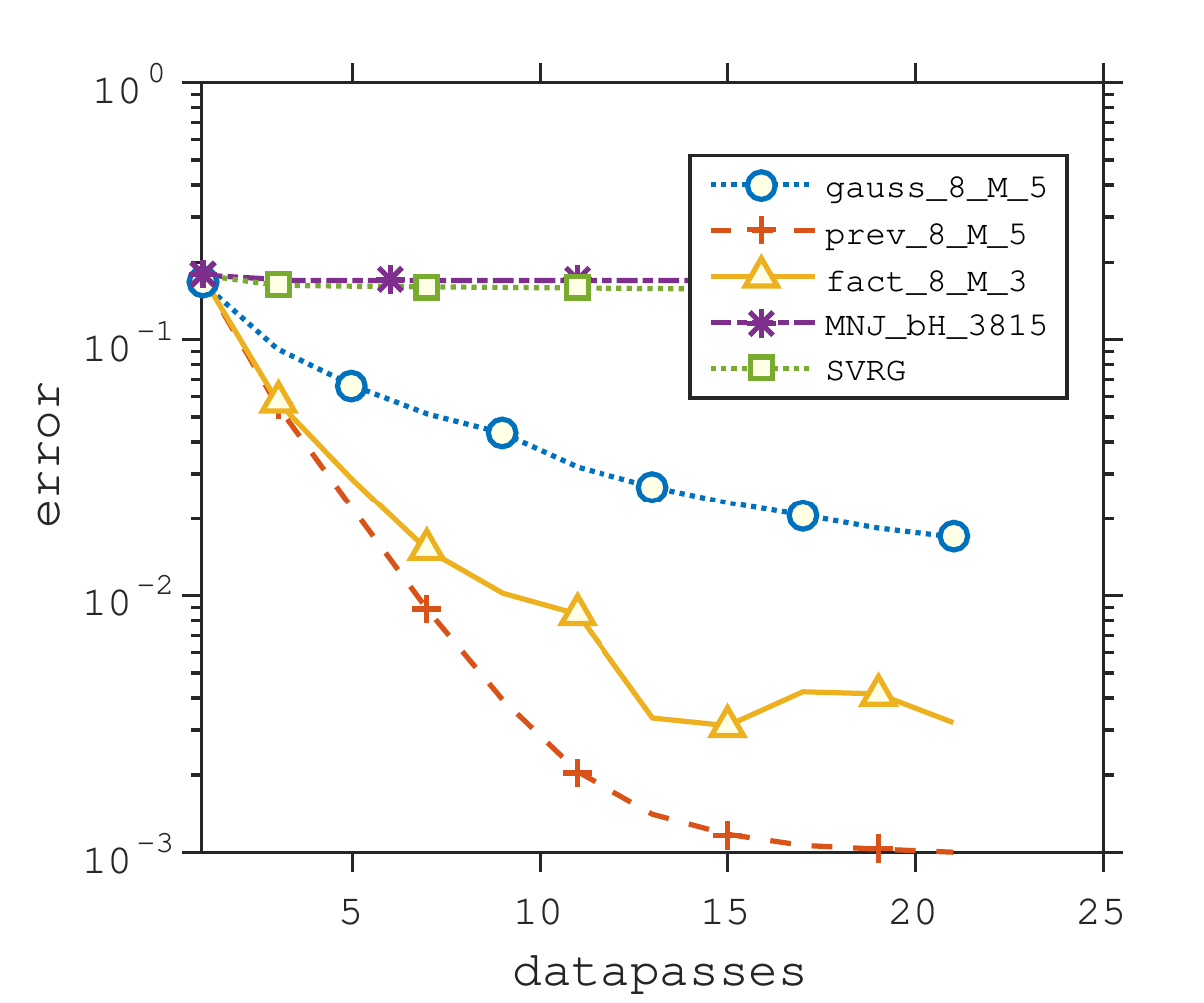}
    }\hfill
\subfigure[\texttt{HIGGS}]{	
\includegraphics[width =  \widshrink\columnwidth, height =\heightshrink\columnwidth]{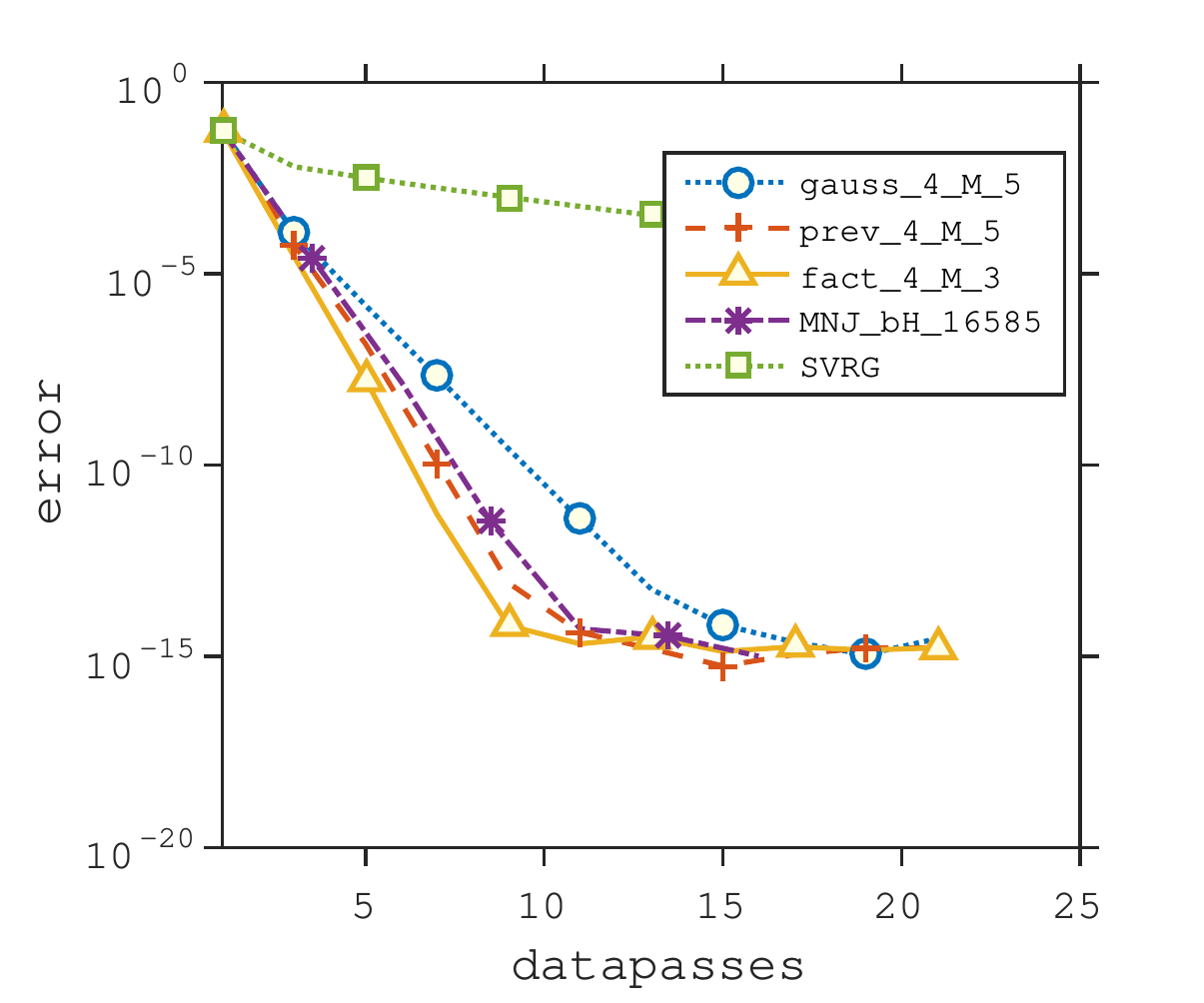}
	} \hfill
\subfigure[\texttt{\texttt{SUSY}}]{
	\includegraphics[width =  \widshrink\columnwidth, height =\heightshrink\columnwidth]{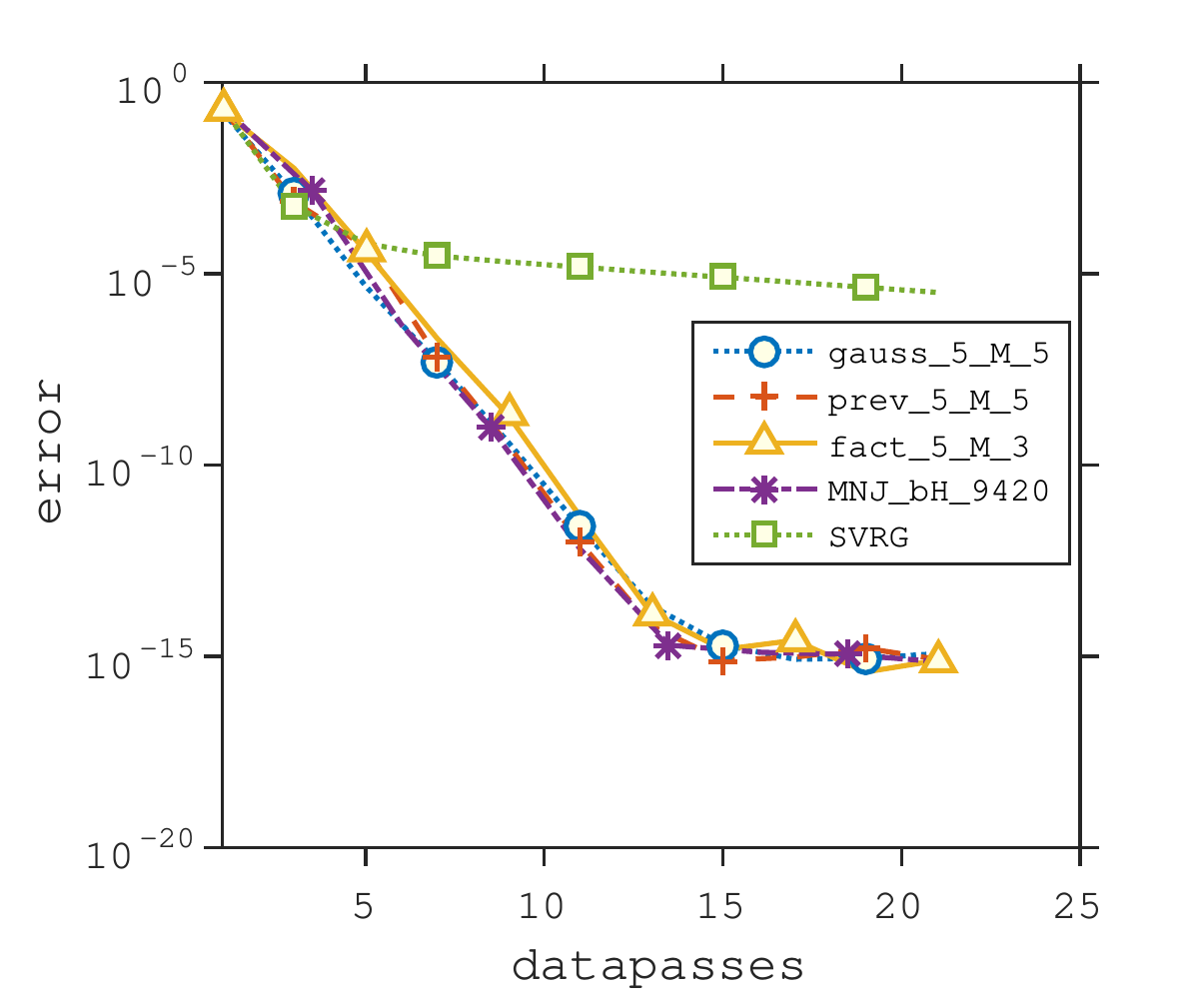}
}\hfill%
    \subfigure[ \texttt{epsilon\_normalized} ]{
    \includegraphics[width =  \widshrink\columnwidth, height =\heightshrink\columnwidth]{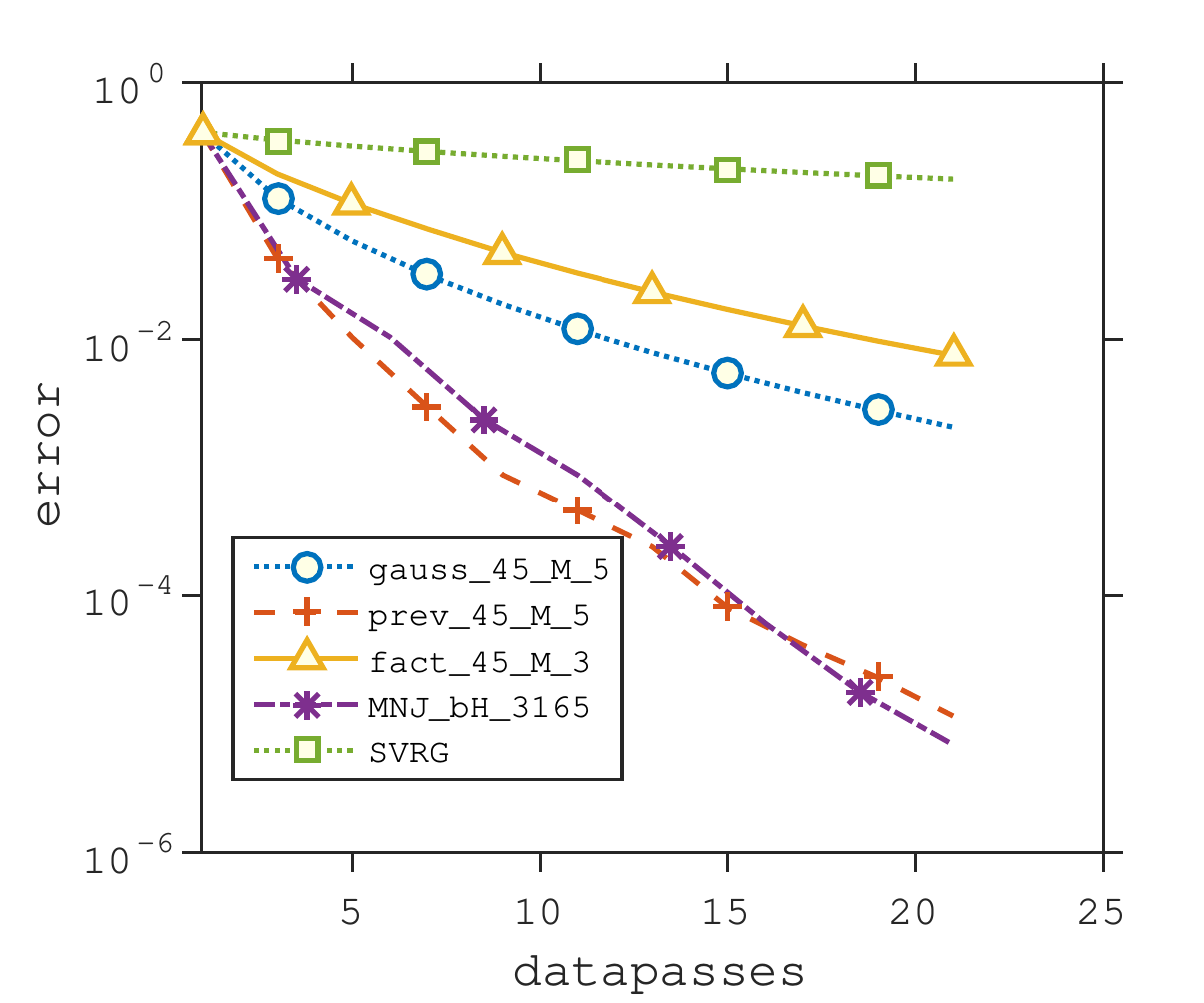}}
    \caption{ (a) \texttt{gissette\_scale} $(d;n) = (5,001; \; 6,000)$ (b) \texttt{covtype-libsvm-binary} $(d;n) = (55; \; 581,012)$ (c) \texttt{HIGGS} $(d;n) = (29; \; 11,000,000)$ (d) \texttt{SUSY} $(d;n) = (19;\; 3,548,466)$  (e)  \texttt{epsilon\_normalized} $(d;n) = (2,001; \;400,000)$ } \label{fig:LIBSVM1}
\end{figure}

\begin{figure}
\centering
\subfigure[\texttt{rcv1-train.binary}]{	
\includegraphics[width =  \widshrink\columnwidth, height =\heightshrink\columnwidth ]{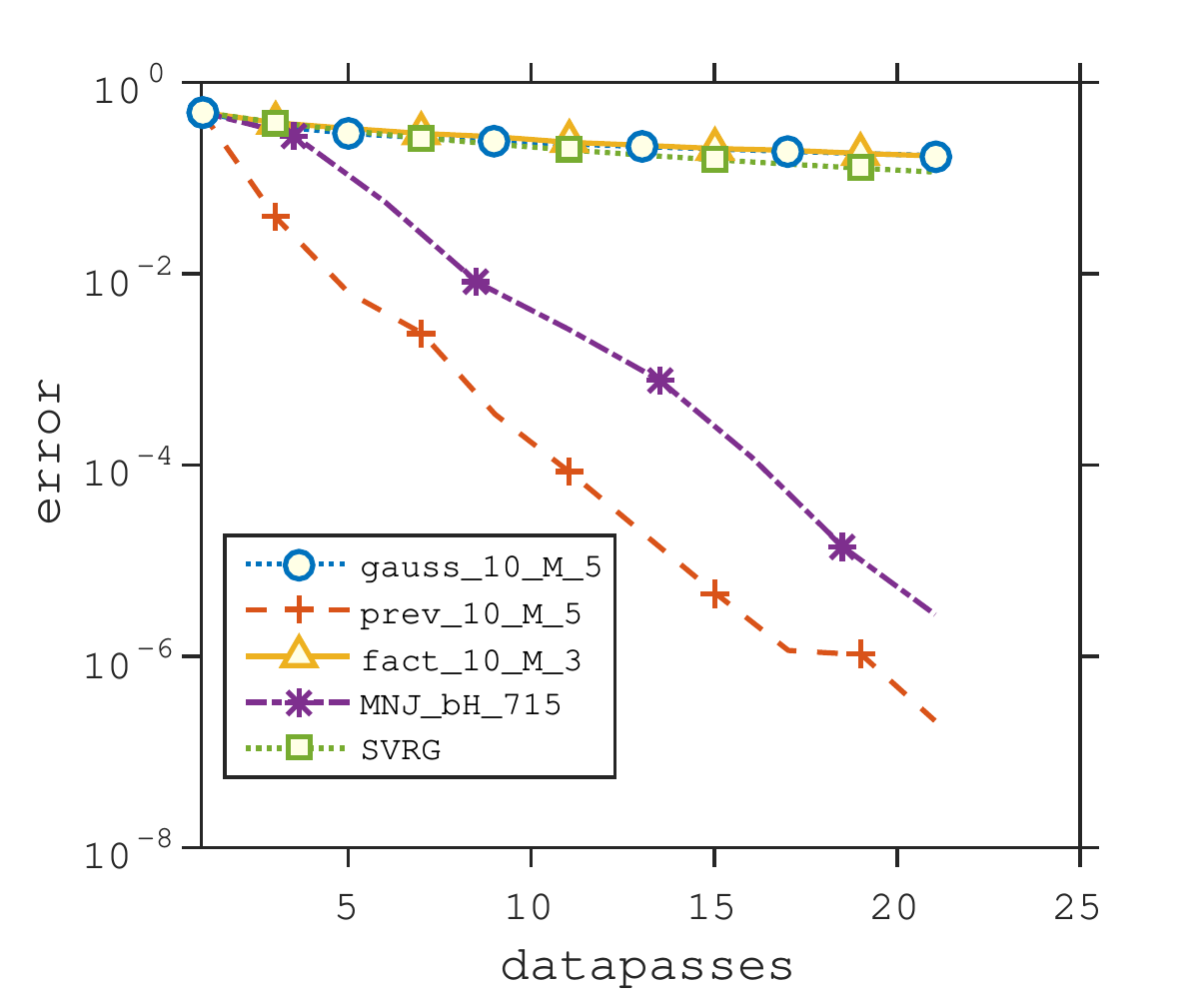}
	} \hfill
\subfigure[\texttt{url-combined}]{
	\includegraphics[width =  \widshrink\columnwidth, height =\heightshrink\columnwidth ]{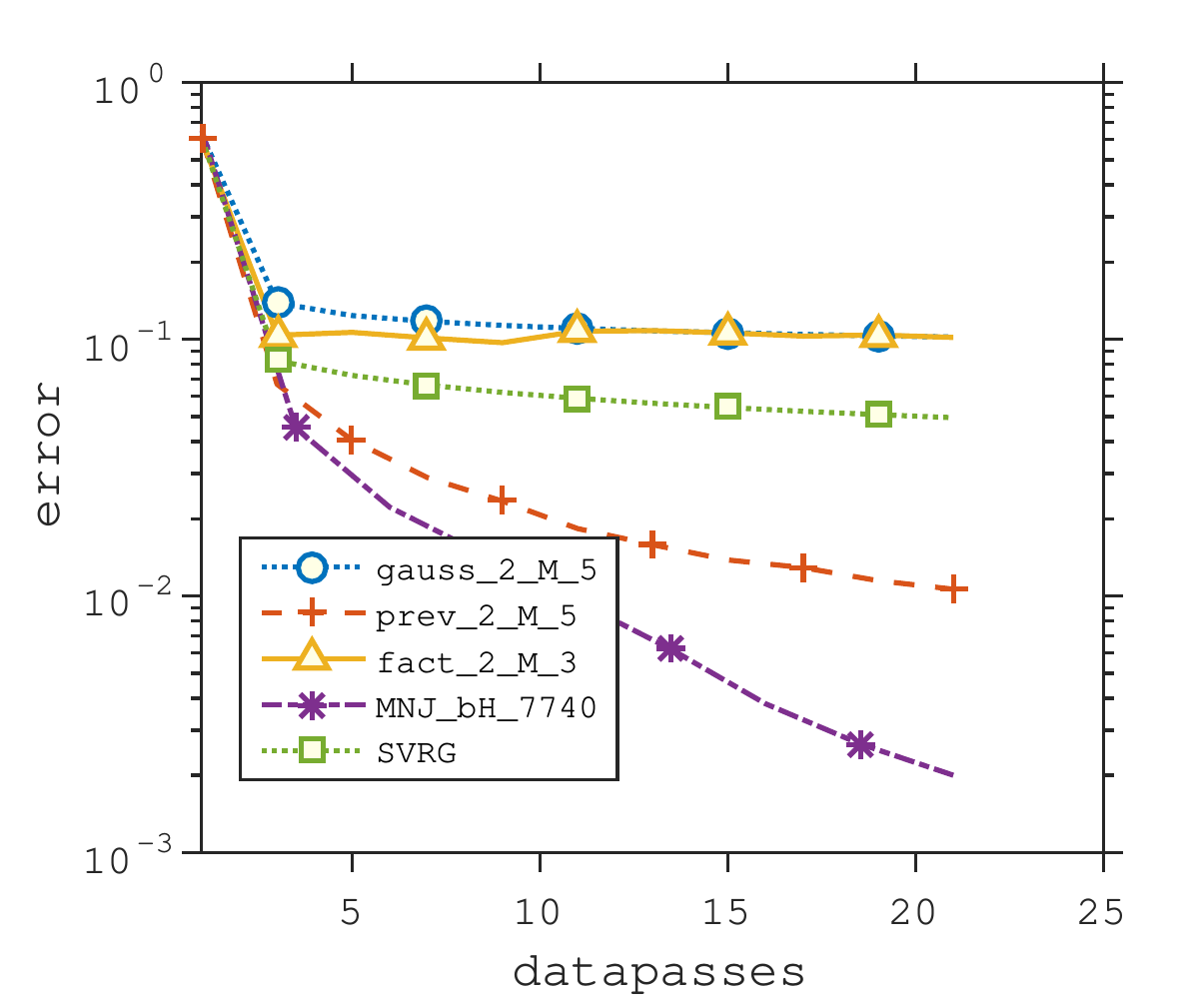}
}%
    \caption{(a) \texttt{rcv1-train.binary} $(d;n) = (47,237; \; 20,242)$ (b) \texttt{url-combined} $(d;n) = (3,231,962; \; 2,396,130)$.} \label{fig:LIBSVM2}
\end{figure}

\subsection{Timed}\label{sec:numtime}

\begin{figure}
 \centering
    \subfigure[ \texttt{gissette\_scale} ]{
    \includegraphics[width =  \widshrink\columnwidth, height =\heightshrink\columnwidth ]{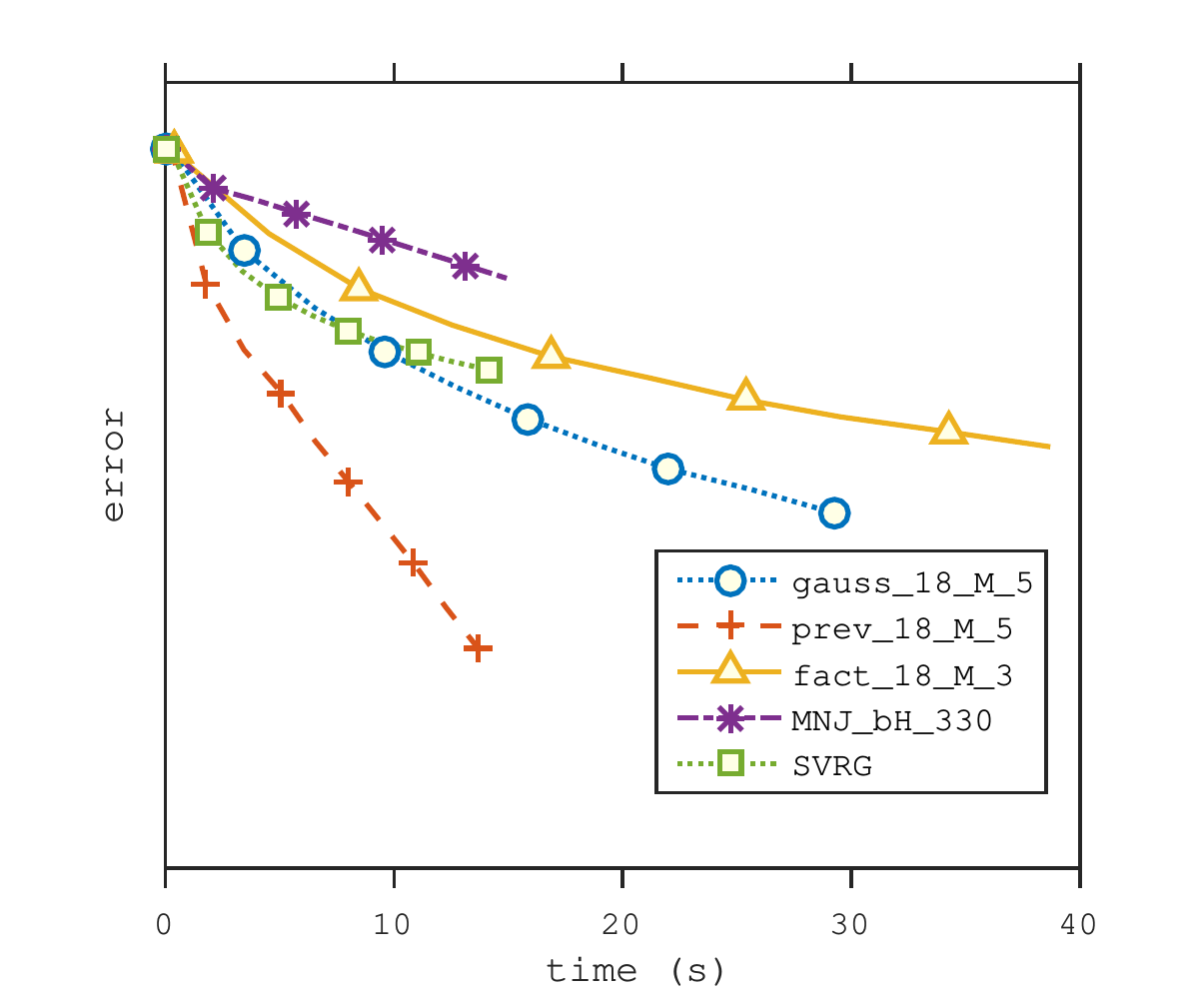}
    }\hfill
    \subfigure[ \texttt{covtype-libsvm-binary} ]{
    \includegraphics[width =  \widshrink\columnwidth, height =\heightshrink\columnwidth ]{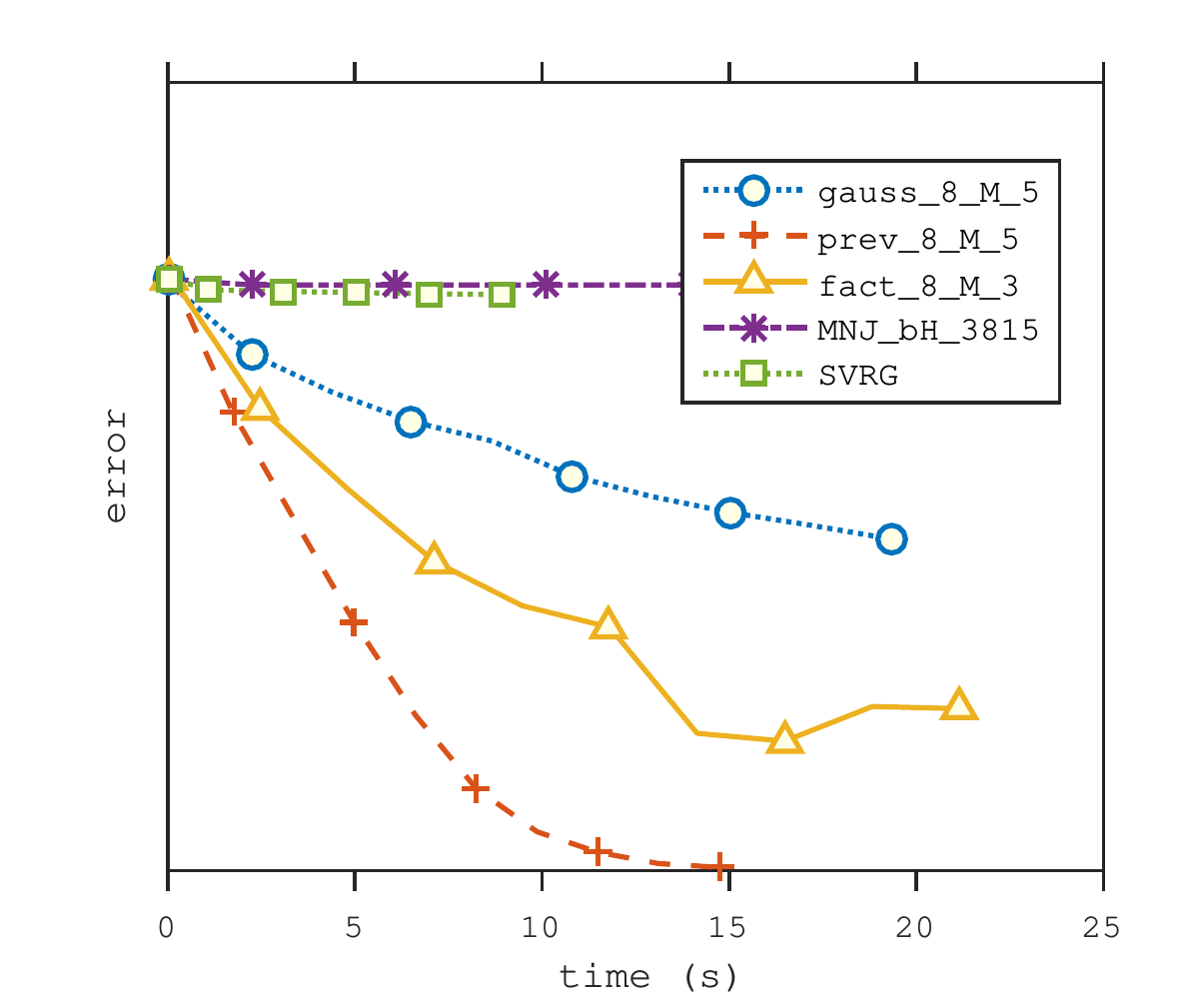}
    }\hfill
\subfigure[\texttt{HIGGS}]{	
\includegraphics[width =  \widshrink\columnwidth, height =\heightshrink\columnwidth]{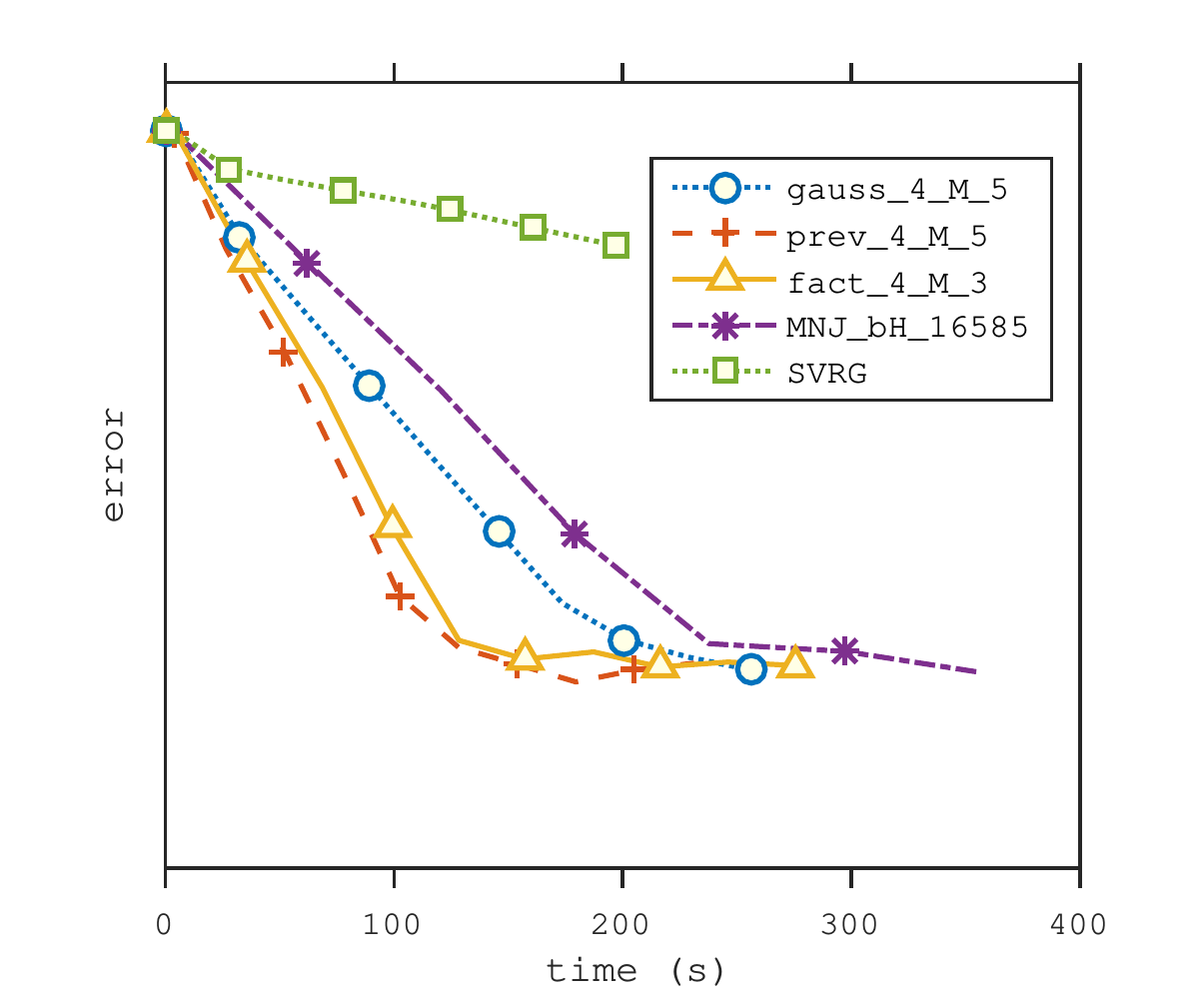}
	} \hfill
\subfigure[\texttt{\texttt{SUSY}}]{
	\includegraphics[width =  \widshrink\columnwidth, height =\heightshrink\columnwidth]{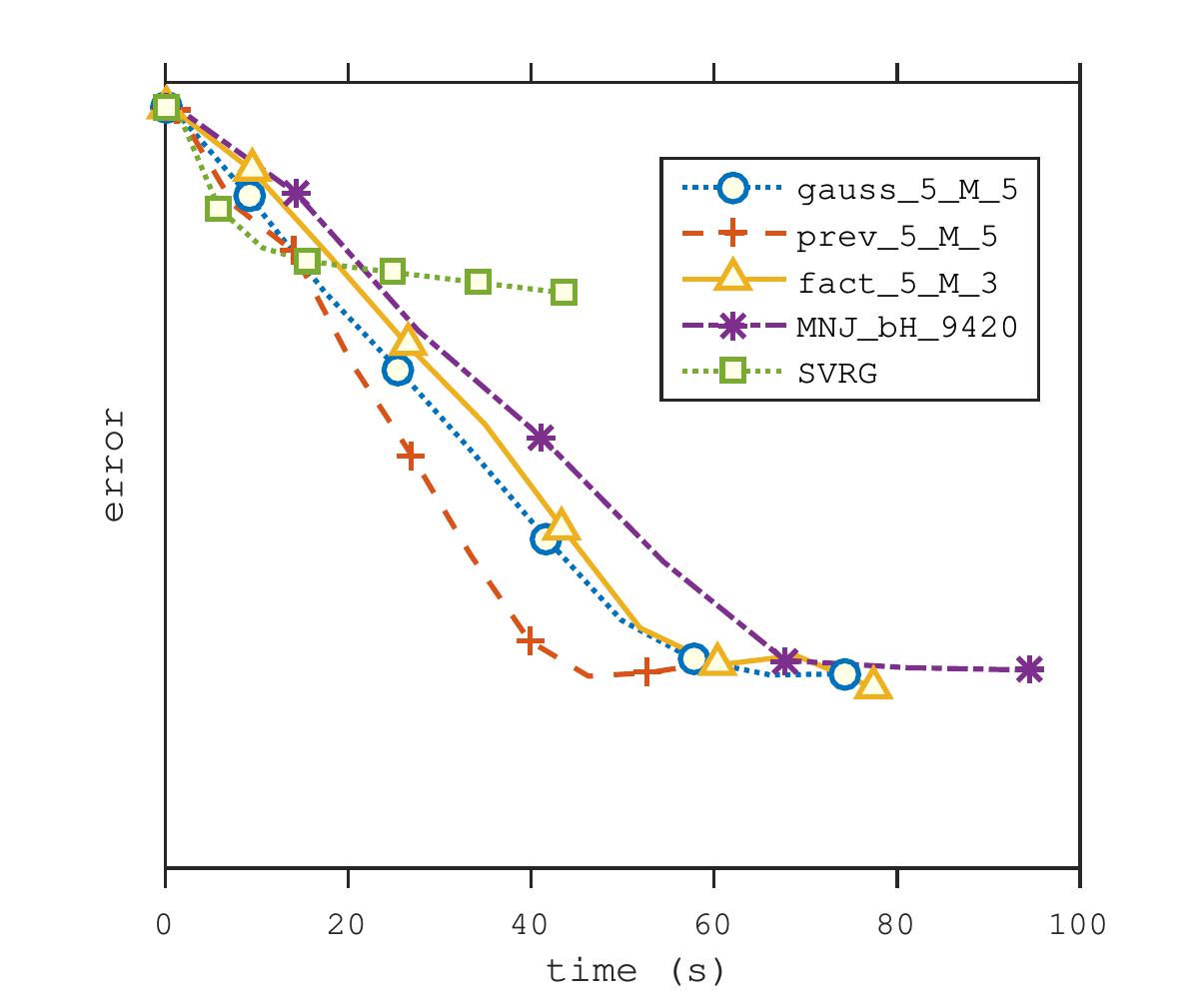}
}\hfill%
    \subfigure[ \texttt{epsilon\_normalized} ]{
    \includegraphics[width =  \widshrink\columnwidth, height =\heightshrink\columnwidth]{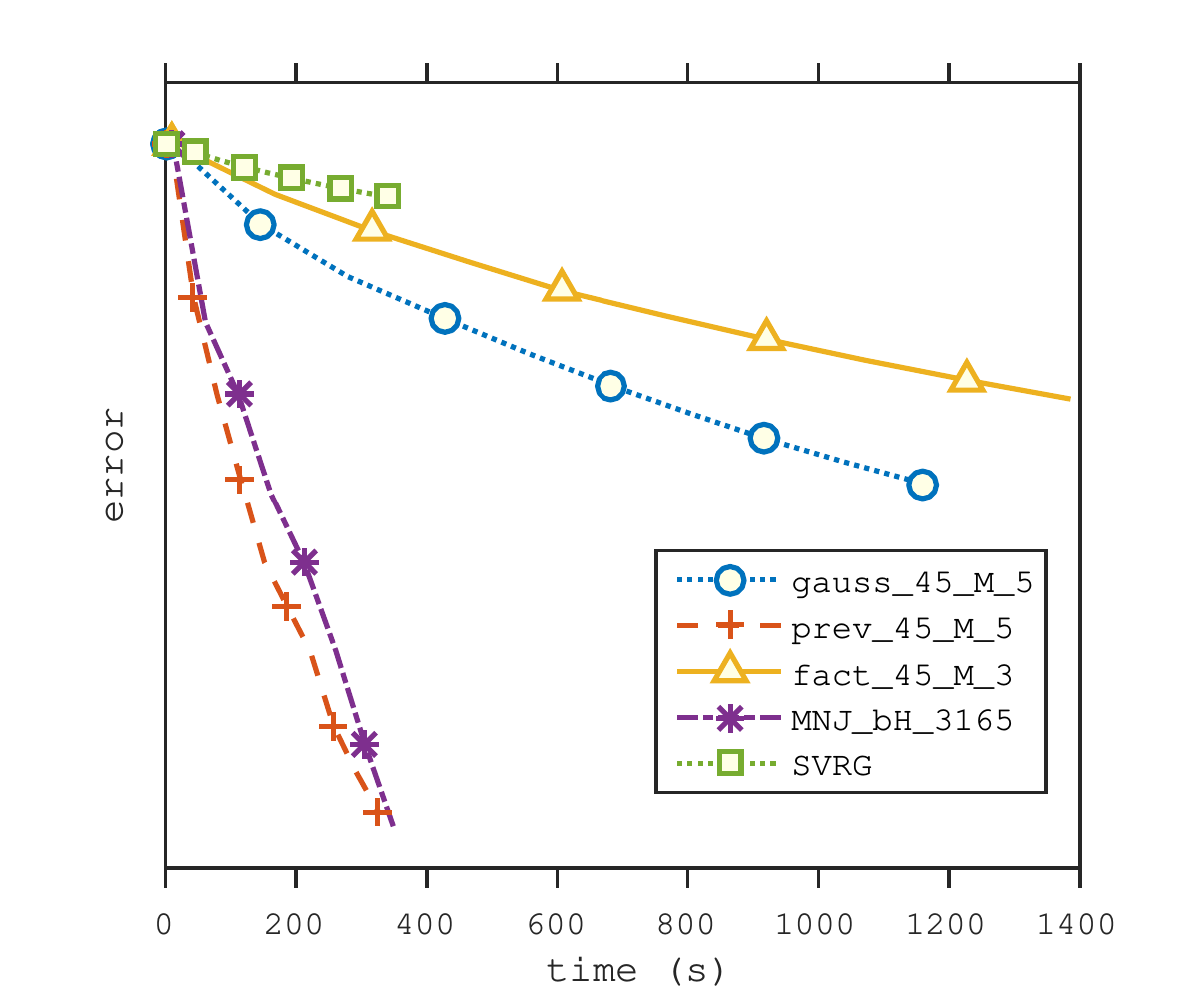}}
    \caption{ (a) \texttt{gissette\_scale} $(d;n) = (5,001; \; 6,000)$ (b) \texttt{covtype-libsvm-binary} $(d;n) = (55; \; 581,012)$ (c) \texttt{HIGGS} $(d;n) = (29; \; 11,000,000)$ (d) \texttt{SUSY} $(d;n) = (19;\; 3,548,466)$  (e)  \texttt{epsilon\_normalized} $(d;n) = (2,001; \;400,000)$ } \label{fig:LIBSVMtime1}
\end{figure}

\begin{figure}
\centering
\subfigure[\texttt{rcv1-train.binary}]{	
\includegraphics[width =  \widshrink\columnwidth, height =\heightshrink\columnwidth ]{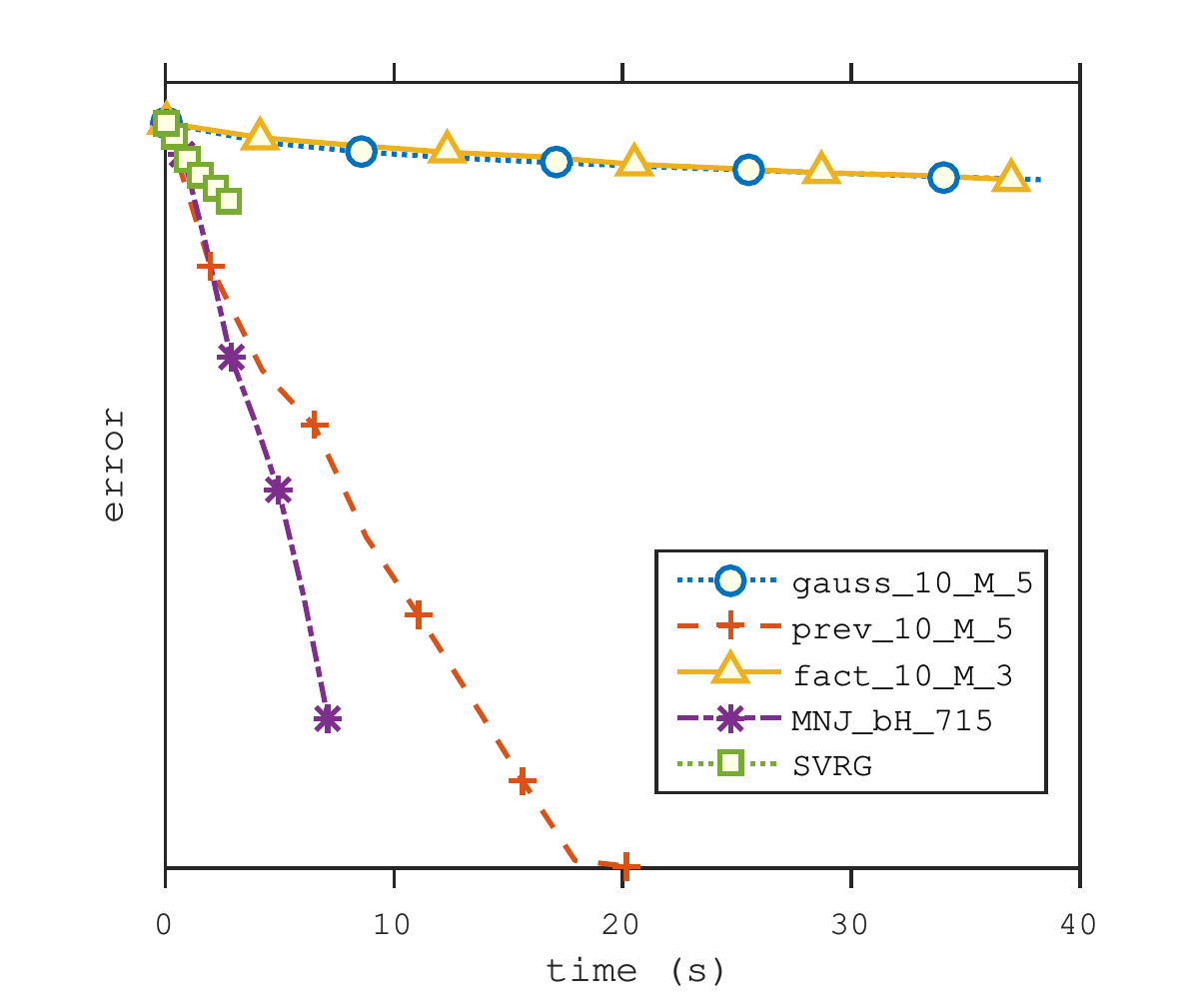}
	} \hfill
\subfigure[\texttt{url-combined}]{
	\includegraphics[width =  \widshrink\columnwidth, height =\heightshrink\columnwidth ]{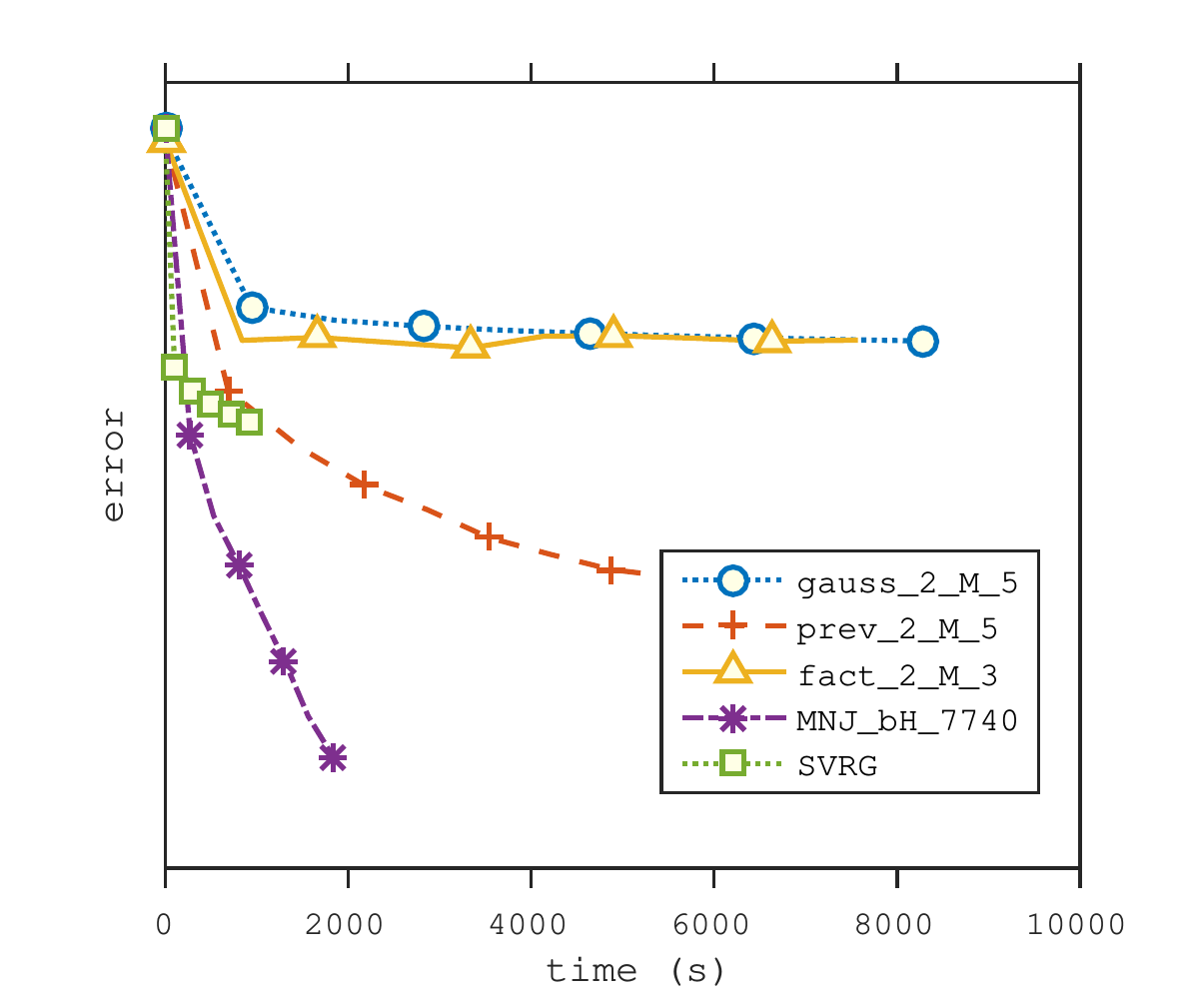}
}%
    \caption{(a) \texttt{rcv1-train.binary} $(d;n) = (47,237; \; 20,242)$ (b) \texttt{url-combined} $(d;n) = (3,231,962; \; 2,396,130)$.} \label{fig:LIBSVMtime2}
\end{figure}

In  Figures~\ref{fig:LIBSVMtime1} and~\ref{fig:LIBSVMtime2}, we
compare the evolution of the error over time for each method.
 While measuring time is implementation and machine dependent, we include these time plots as to 
  provide 
 further insight into the methods performance. Note that we did not use any sophisticated implementation tricks such as ``lazy'' gradient updates~\cite{S2GD}, but instead implemented each method as originally designed so that the methods can be compared on a equal footing.

The results in these tests corroborate with our conclusions in Section~\ref{sec:numdatapass}

\section{Extensions} This work opens up many new research avenues. For instance, sketching techniques are increasingly successful tools in large scale machine learning, numerical linear algebra, optimization and computer science. Therefore, one could employ a number of new sketching methods in the block BFGS method, such as the Walsh-Hadamard matrix~\cite{Pilanci2014,Lu2013}. Using new sophisticated sketching methods, combined with the block BFGS update, could result in even more efficient and accurate estimates of the underlying curvature.


Viewed in terms of sketching, the \texttt{MNJ} method~\cite{Moritz2015} and our Stochastic Block BFGS method with \texttt{prev} sketching  follow opposite strategies.
While the \texttt{MNJ} method sketches the Hessian matrix with a single vector $v$, where $v$ is the average over a combination of previous search directions (a very coarse approximation), our  \texttt{prev} variant
 uses all previous search directions to form the sketching matrix $D_t$ (a finer approximation).
Deciding between these two extremes or how to combine them could be done adaptively by  examining the curvature matrix $D_t \nabla f (x_t) D_t.$ When the previous search directions are almost collinear,  $D_t \nabla f (x_t) D_t$ becomes ill-conditioned. In this case one should form matrix $D_t$ with less columns using fewer or a coarser combination of previous search directions, while if $D_t \nabla f (x_t) D_t$ is well-conditioned, one could use  more or a finer combination of previous search directions.

While in this work  we have for simplicity focused on utilizing our metric learning techniques in conjunction with SVRG, they can be used with other optimization algorithms as well, including SGD~\cite{RobbinsMonro:1951}, SDCA~\cite{SDCA} and more.

\section{Proof of Lemma~\ref{lem:pretheo}}
Let $h_S(w) = f_S(w)-f_S(w_*) -\nabla f_S(w_*)^T(w-w_*)$. Note that $h_S(w)$ achieves its minimum at $w_*$ and $h_S(w_*)=0.$ Furthermore,  $\nabla^2 h_S(w) \preceq \Lambda I$ from Assumption~\ref{ass:strongsmooth}. Consequently,  for every $b \in \R^d$ we have
\begin{align*}
0&= h_S(w_*) \leq h_S(w+b) \leq  h_S(w) +\nabla h_S(w)^Tb + \frac{\Lambda }{2} \norm{b}_2^2.
\end{align*}
Minimizing the right hand side of the above in $b$ gives
\begin{align*}
0&=h_S(w_*) \leq  h_S(w) -  \frac{1}{2\Lambda} \norm{\nabla h_S(w)}_{2}^2.
\end{align*}
Re-arranging the above and switching back to $f$, we have
\begin{align*}
&\norm{\nabla f_S(w)-\nabla f_S(w_*)}_{2}^2 \leq  2\Lambda\left(f_S(w)-f_S(w_*) -\nabla f_S(w_*)^T(w-w_*)\right).
\end{align*}
Recalling the notation $\df(x) \eqdef f(x) - f( w_*)$ and taking expectation with respect to $S$ gives
\begin{equation} \label{eq:xxxx}\E{\norm{\nabla f_S(w)-\nabla f_S(w_*)}_{2}^2  } \leq
2\Lambda \df(w).
\end{equation}
Now we apply the above to obtain
\begin{eqnarray*}
\E{\norm{g}_{2}^2 } & \leq &  2\E{\norm{\nabla f_S(x)-\nabla f_S(w_*)}_{2}^2} + 2\E{\norm{ \nabla f_S(w) - \nabla f_S(w_*) -\mu}_{2}^2}\\
& \leq & 2\E{\norm{\nabla f_S(x)-\nabla f_S(w_*)}_{2}^2}+ 2\E{\norm{ \nabla f_S(w) - \nabla f_S(w_*)}_{2}^2} - 2\norm{\nabla f(w)}_{2}^2.
\end{eqnarray*}
where we used $\norm{a+b}_2^2 \leq 2\norm{a}_2^2+2\norm{b}_2^2$ in the first inequality and $\mu =\E{\nabla f_S(w_{k}) - \nabla f_S(w_*)}$. Finally,
\begin{eqnarray*}
\E{\norm{g}_{2}^2  }
&\overset{\eqref{eq:xxxx}}{\leq} &4\Lambda\left( \df(x)+ \df(w)\right) - 2\norm{\nabla f(w)}_{2}^2\\
&\leq   &4\Lambda \df(x) + 4(\Lambda-\lambda)\df(w).
\end{eqnarray*}
In the last inequality we used the fact that strongly convex functions satisfy  $\|\nabla f(x)\|_2^2~\geq~2\lambda~\df(x)$ for all $x\in \R^d.$

\section{Proof of Lemma~\ref{lem:Hspectra}}
To simplify notation, we define $ G \eqdef \nabla^2 f_T (x_t)$,
$\Delta =\Delta_t $,
$Y = Y_t $, $H= H_{t-1}$, $H^+ = H_t$, $B=H^{-1}$, $B^+ = (H^+)^{-1}$ and
$V=Y \Delta D^T$.
Thus, the block BFGS update can be written as
$$H^+ = H - V^T H - H V +V^T H V +D \Delta D^T.$$

Proposition 2.2 in~\cite{Gower2014c} proves that so long as $G$ and $H$ (and hence $B$) are positive definite and $D$ has full rank, then $H^{+}$ is positive definite and non-singular, and consequently, $B^+$ is well defined and positive definite.
Using the Sherman-Morrison-Woodbury identity the update formula for $B^+$, as shown in the Appendix in~\cite{Gower2014c}, is given by
\begin{equation}\label{B}
B^+ = B + Y \Delta Y\t  -
B D (D\t B D)^{-1} D\t B.
\end{equation}
We will now bound $\lambda_{\max}(H^+)=\|H^+\|_2$ from above and $\lambda_{\min}(H^+)=1/\|B^+\|_2$ from below.

  Let $C = B D ( D\t B D )^{-1} D\t B $. Then since $C \succeq 0 $, $B-C \preceq B$ and hence, $\|B-C\|_2 \leq \|B\|$ and
$$\|B^+\|_2 \leq \|B\|_2 + \| Y \Delta Y\t \|_2. $$
Now, letting $G^\half$ and $G^{-\half}$ denote the unique square root of $G$ and its inverse, and defining $U=G^\half D$, we have
$$ D \Delta D\t
= G^{-\half} U(U\t U)^{-1}U\t G^{-\half} =
G^{-\half} P G^{-\half} ,$$
where $P = U(U\t U)^{-1}U\t$ is an orthogonal projection matrix. Moreover, it is easy to see that
$$ Y \Delta Y\t
= G^\half P G^\half
\ \ {\rm and} \ \ V =Y \Delta D\t =  G^\half P G^{-\half}.$$
Since $\|MN\|_2 \leq \|M\|_2 \|N\|_2$ and $\|P\|_2 = 1$, we have
$\|D \Delta D\t \|_2 \leq \|G^{-1}\|_2$, $\|Y \Delta Y\t \|_2 \leq \| G \|_2$  and $ \|Y \Delta D\t \|_2 \leq \|G^{-\half}\|_2\|G^\half\|_2$.

Hence,
\begin{equation}\label{eq:Bplusbnd}
\|B^+\|_2 \leq \|B\|_2 + \| G \|_2 \overset{\eqref{eq:hessbound}}{\leq} \|B\|_2 + \Lambda.
\end{equation}
Furthermore,
\begin{eqnarray}
\|H^+\|_2  &\leq & \|H\|_2 + 2 \|H\|_2 \|G^{-\half}\|_2 \|G^\half \|_2 \nonumber
 + \|H\|_2 \|G^{-1}\|_2 \| G \|_2 + \|G^{-1}\|_2 \nonumber \\
& \leq & (1 + 2\sqrt{\kappa} + \kappa)\|H\|_2 + \frac{1}{\lambda} \\
&=& \alpha \|H\|_2 + \frac{1}{\lambda}, \label{eq:normHbnd}
\end{eqnarray}
where $\kappa = \Lambda/\lambda$ and $\alpha = (1+\sqrt{\kappa})^2$.

Since we use a memory of $M$ block triples $(D_i, Y_i, \Delta_i)$, and the metric matrix $H_t$ is the result of applying, at most, $M$ block updates  BFGS~\eqref{eq:BBFGS} to $H_0$,
we have that
\begin{equation} \label{eq:Trbnbt}
\lambda_{\max}(B_t) = \|B_{t}\|
\overset{\eqref{eq:Bplusbnd}}{\leq} \|B_{t-M}\| + M  \Lambda,
\end{equation}
and hence that
\[ \gamma = \lambda_{\min}(H_t) \overset{\eqref{eq:Trbnbt}}{\geq} \frac{1}{ \|B_{t-M}\| + M \Lambda}.\]
Finally, since $\alpha = (1+\sqrt{\kappa})^2$, we have
\begin{eqnarray*}
\Gamma &= & \lambda_{\max}(H_t) \\
& = & \|H_t\|\\
 & \overset{\eqref{eq:normHbnd}}{ \leq} & \alpha^M \|H_{t-M}\| + \frac{1}{\lambda}\sum_{i=0}^{M-1}\alpha^{i} \\
&= & \alpha^M \|H_{t-M}\| + \frac{1}{\lambda}\frac{\alpha^M - 1}{\alpha - 1} \\
&\leq & (1+\sqrt{\kappa})^{2M} \left(\|H_{t-M} \| + \frac{1}{\lambda (2\sqrt{\kappa} + \kappa)}\right).
\end{eqnarray*}
The bound~\eqref{eq:Gammakappbnd} now follows by observing that $H_{t-M} =I.$

{ \small
\printbibliography
}
\end{document}